\documentclass[11pt,reqno]{amsart}
\usepackage[utf8]{inputenc}
\usepackage{amssymb}
\usepackage{amsmath} %\allowdisplaybreaks
\usepackage{amsthm}
\usepackage{amsfonts}
\usepackage[textsize=small]{todonotes} % put 'disable' in the brackets to remove all todonotes
\usepackage{hyperref}
\usepackage{a4wide}
\usepackage{dsfont} % for \mathds{1}
\usepackage{subcaption}
\hypersetup{
    colorlinks,
    citecolor=red,
    filecolor=black,
    linkcolor=blue,
    urlcolor=black
}

\newtheorem{theorem}{Theorem}[section]
\newtheorem{lemma}[theorem]{Lemma}
\newtheorem{proposition}[theorem]{Proposition}

\theoremstyle{definition} \newtheorem{definition}[theorem]{Definition}
\theoremstyle{definition} \newtheorem{remark}[theorem]{Remark}
\theoremstyle{definition} \newtheorem*{remark*}{Remark}

\newtheorem{corollary}[theorem]{Corollary}

\renewcommand{\epsilon}{\varepsilon}
\renewcommand{\tilde}{\widetilde}

\DeclareMathSymbol{\shortminus}{\mathbin}{AMSa}{"39}

\newcommand{\rd}{{\textup{d}}}

\newcommand{\R}{{\mathbb{R}}}
\newcommand{\C}{{\mathbb{C}}}
\newcommand{\Z}{{\mathbb{Z}}}
\newcommand{\N}{{\mathbb{N}}}

\newcommand{\cC}{{\mathcal{C}}}

\newcommand{\g}{{\tilde{g}}}

\title{Families of periodic delay orbits}
\author{Peter Albers, Philipp Aretz, Irene Seifert}
\address{Mathematisches Institut, Universit\"at Heidelberg,
Im Neuenheimer Feld 205, 69120 Heidelberg, Germany}
\email{palbers@mathi.uni-heidelberg.de}
\email{philipp.aretz@stud.uni-heidelberg.de}
\email{iseifert@mathi.uni-heidelberg.de}
\date{\today}

\begin{document}

\begin{abstract}
We construct and analyze families of periodic delay orbits for a class of delay differential equations in two dimensions depending on two real-valued functions. These families are parametrized by the delay parameter. It is possible to represent the dependency of these periodic delay orbits on the delay parameter by a curve in the plane, without loss of information. It turns out that the singularities of these curves necessarily are cusps in the non-degenerate case. After discussing degenerate situations in general, we explain how to glue different families of periodic delay orbits at degeneracies in the delay parameter. 

%We construct a wide class of examples for families of periodic delay equations parameterized by delay. Analyzing the examples, we find out more about the structure of these families. We analyze, for instance, what happens when the parameterization by delay attains zero speed, and the case of two smooth families crossing in a degenerate periodic orbit. The examples can be used as a source of inspiration for further research.

%    \vspace{0.2cm}
%    
%    \textbf{Keywords:}
%    delay differential equations,
%    periodic orbits
\end{abstract}
    
\maketitle

%\tableofcontents

\section{Introduction}

Many dynamical systems are modelled on differential equations of the type
\begin{align}
    \dot{x}(t) = X_t\big( x(t) \big), \label{eq:ODE-intro}
\end{align}
where we interpret $t\in\R$ as a time parameter. Here, $X$ is a smooth time-dependent vector field on $\R^n$ (or, more generally, on a manifold) and $x:\R\to\R^n$ is a differentiable function with derivative $\dot{x}$. A solution $x$ of \eqref{eq:ODE-intro} (which is necessarily smooth by a bootstrapping argument) is called an \textit{orbit} of the vector field $X$. \textit{Periodic} orbits, i.e.\ those that satisfy $x(t+T)=x(t)$ for some $T\in\R$ and all $t\in\R$, are of particular interest.

In many situations, however, it is natural to consider a delay differential equation of the type
\begin{align}
    \dot{x}(t) = X_t\big( x(t-\tau) \big), \label{eq:DDE-intro}
\end{align}
where we have an additional parameter $\tau\in\R$ interpreted as the delay between the state of the system and the reaction to it. We call a solution of \eqref{eq:DDE-intro} a $\tau$-\textit{delay orbits} and are particularly interested in \textit{periodic delay orbits}.

Although looking similar to ordinary differential equations, delay differential equations are significantly more difficult to treat. For instance, while an initial condition for \eqref{eq:ODE-intro} is simply a point in $\R^n$, an initial condition for \eqref{eq:DDE-intro} consists of a whole function on an interval of length $\tau$, called the \textit{initial history}. This means that, even for formulating an initial value problem, we have to work with infinite-dimensional spaces. Moreover, delay differential equations do not admit a flow, instead one considers a \textit{semi-flow} or a \textit{forward evolutionary system} on a suitable function space. In particular, finding periodic orbits -- for instance using fixed point theorems for the flow map -- is much harder for delay differential equations. Moreover, the shift map $x(\cdot)\mapsto x(\cdot - \tau)$ is not smooth between suitable function spaces (cf.~\cite{FrauenfelderWeber18}), so that it is not straightforward to consider \eqref{eq:DDE-intro} as a perturbation resp.~deformation of \eqref{eq:ODE-intro}.

A good overview of the theory of delay differential equations can be found in the book \cite{Diekmann95}. 
The article \cite{AlbersSeifert22} is written more from a symplectic dynamics viewpoint and focuses on periodic delay orbits. It uses sc-smoothness, a part of polyfold theory \cite{HoferWysockiZehnder21}, and, in particular, the polyfold implicit function theorem to make the idea precise that a delay differential equation \eqref{eq:DDE-intro} can indeed be treated as a perturbation of the ordinary differential equation. This leads to the following persistence result for non-degenerate periodic orbits.

\begin{theorem}[{\cite[Theorem 1.1]{AlbersSeifert22}}]
\label{thm:mytheorem}
Let $x_0$ be a non-degenerate 1-periodic orbit of a smooth time-dependent vector field $X: S^1\times\R^n\to\R^n$. Then there exists $\tau_0 > 0$ such that for every delay $\tau$ with $|\tau|\leq \tau_0$ there exists a (locally unique) smooth 1-periodic solution $x_{\tau}$ of the delay equation $\dot{x}(t)=X_t\big( x(t-\tau)\big)$. Moreover, the parameterization $\tau\mapsto x_{\tau}$ is smooth.
\end{theorem}

The $\tau$-delay orbit $x_{\tau}$ from Theorem~\ref{thm:mytheorem} lies near the original orbit $x_0$ in the $\cC^{\infty}$-topology, i.e.~$x_0, x_{\tau}: S^1 \to \R^n$ are close as maps with all their derivatives. Also, local uniqueness is to be understood in $\cC^{\infty}$-topology.
For the definition of non-degeneracy, see Section \ref{sec:periodic-orbits} below.
The statement from Theorem~\ref{thm:mytheorem} continues to hold for much more general delay equations than \eqref{eq:DDE-intro}, see \cite{AlbersSeifert22} and \cite{Seifert22}.

Theorem \ref{thm:mytheorem} tells us that every non-degenerate 1-periodic orbit $x_0$ is part of a smooth 1-parameter family $x_\tau$ of delay orbits, and that periodic delay orbits $x_\tau$ for small delay $\tau$ look similar to classical periodic orbits. However, Theorem \ref{thm:mytheorem} is not constructive. Neither does it tell us what the delay orbits look like, nor does it give any quantification of how big the delay parameter $\tau$ can be made for a specific vector field $X$ and a specific orbit $x_0$ for its assertion to hold. All this is due to the use of the polyfold implicit function theorem in its proof.

\medskip

One purpose of this article is to analyze, in detail, examples where Theorem \ref{thm:mytheorem} applies and try to find the family of delay orbits explicitly -- and also look at examples where Theorem \ref{thm:mytheorem} does not apply, to see how its assertion then fails. More specifically, we present a large class of examples in $\R^2$ which, in special cases, first appeared in the doctoral thesis \cite{Seifert22} of the third author. 

Another purpose is to simply give a fairly large class of examples of delay orbits in $\R^2$ with explicit periodic delay orbits. It turns out that we can represent these periodic delay orbits by a single value, i.e.~a point in $\R^2$, even though one would expect an entire initial history. This makes it possible to study rather directly the dependence of the periodic delay orbits on the delay parameter. This leads to curves in $\R^2$ with the surprising property that their singularities are necessarily cusps, in the non-degenerate case. In the degenerate case the local uniqueness in Theorem \ref{thm:mytheorem} may fail as we show in examples. Taking advantage of this failure, we can glue various local families of periodic delay orbits. This is discussed in general and then implemented in specific examples. Throughout this article, we illustrate our results by numerical simulations.

\medskip

\textbf{Acknowledgements: } We thank Lucas Dahinden for the suggestion of the very first idea for these examples. We acknowledge funding by the Deutsche Forschungsgemeinschaft (DFG, German Research Foundation) through Germany’s Excellence Strategy EXC-2181/1 - 390900948 (the Heidelberg STRUCTURES Excellence Cluster), the Transregional Collaborative Research Center CRC/TRR 191 (281071066) and the Research Training Group RTG 2229 (281869850).

\section{The main vector field and its 1-periodic orbits}
\label{sec:periodic-orbits}

We set $S^1:=\R/\Z$. Throughout this article, we fix smooth functions 
\begin{align*}
f: S^1\to\R \qquad\text{and}\qquad g: \R_{\geq 0}\to \R \quad\text{with}\quad g(0)=0.
\end{align*} 
On $\R^2\cong\C$ we then define a time-dependent vector field $X$ by
\begin{align}
    &X: S^1 \times \R^2 \longrightarrow \R^2 \nonumber \\
    &X_t(z)=g(\Vert z\Vert) \cdot \frac{z}{\Vert z\Vert} + f\left(\frac{\arg(z)}{2 \pi} - t\right) \cdot 2\pi i   z
    \label{eq:vector-field}
\end{align}
Since $g(0)=0$ the vector field extends continuously to $0\in\R^2$ by setting by $X_t(0)=0$. We call $g(\Vert z\Vert)  \frac{z}{\Vert z\Vert}$ the  \textit{radial} part and $f\left(\frac{\arg(z)}{2 \pi} - t\right) 2\pi i  z$ the \textit{angular} part of $X_t(z)$. The radial part is time-independent, determined by $g$ and radially symmetric, while the angular part is time-dependent, determined by $f$ and equivariant with respect to multiplication by real scalars.
\begin{figure}[h!]
\captionsetup[subfigure]{justification=centering}
  \begin{subfigure}[t]{.32\linewidth}
    \centering
    \includegraphics[width=\linewidth]{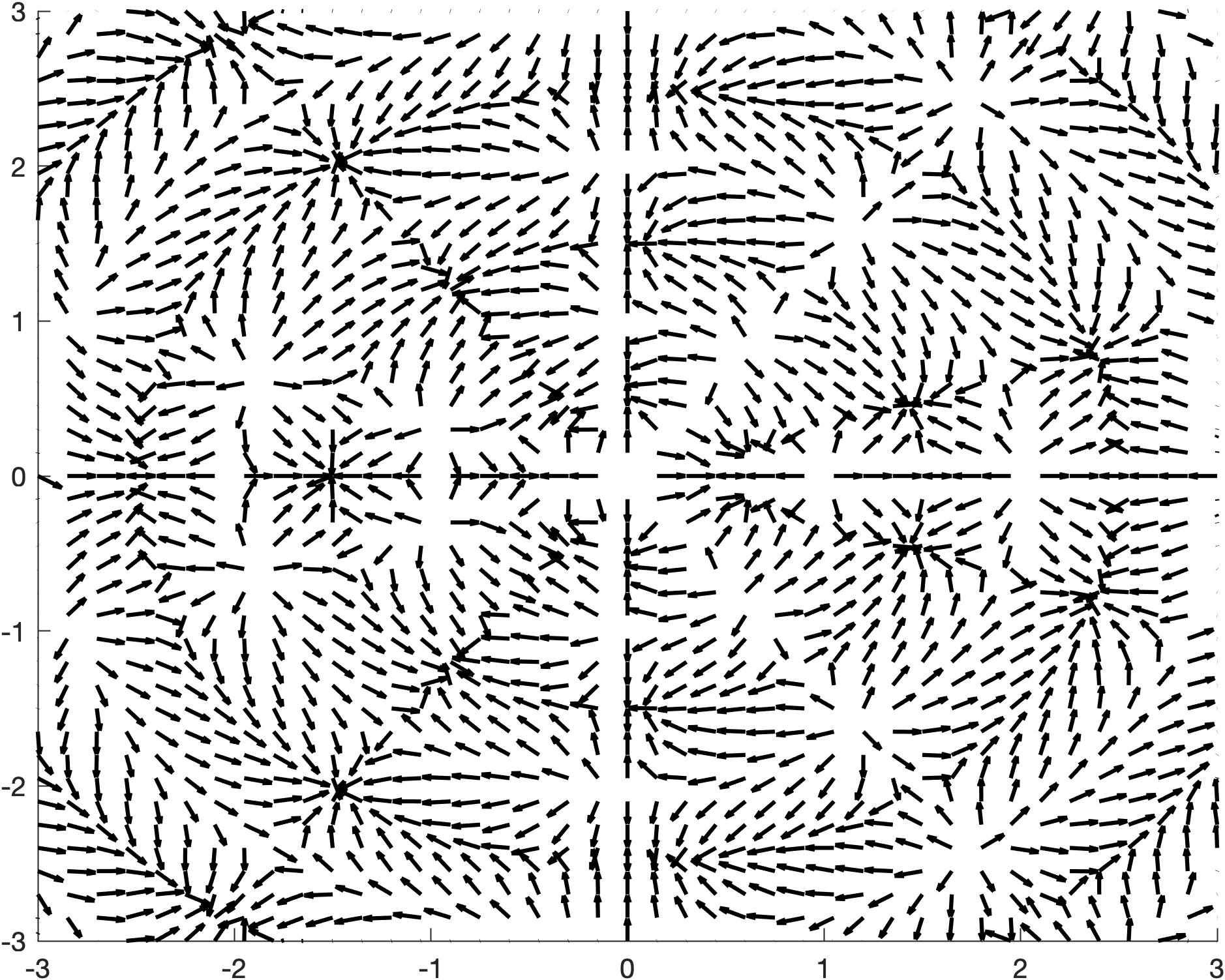}
    \caption{$t=0$ }
    \label{fig:example-vfield-a}
  \end{subfigure}
  \begin{subfigure}[t]{.32\linewidth}
    \centering
    \includegraphics[width=\linewidth]{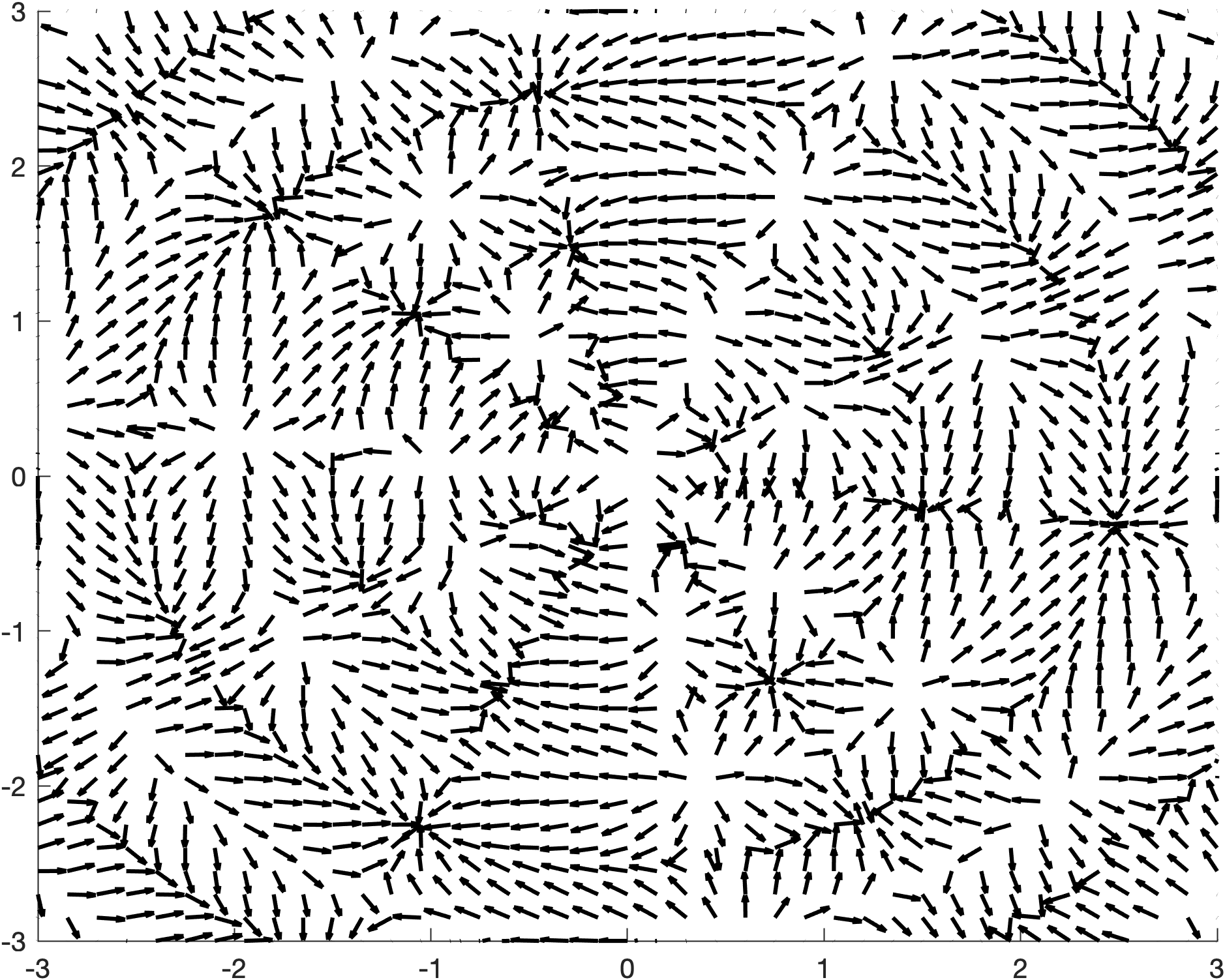}
    \caption{$t=0.33$ }
    \label{fig:example-vfield-b}
  \end{subfigure}
    \begin{subfigure}[t]{0.32\linewidth}
        \centering
        \includegraphics[width=\linewidth]{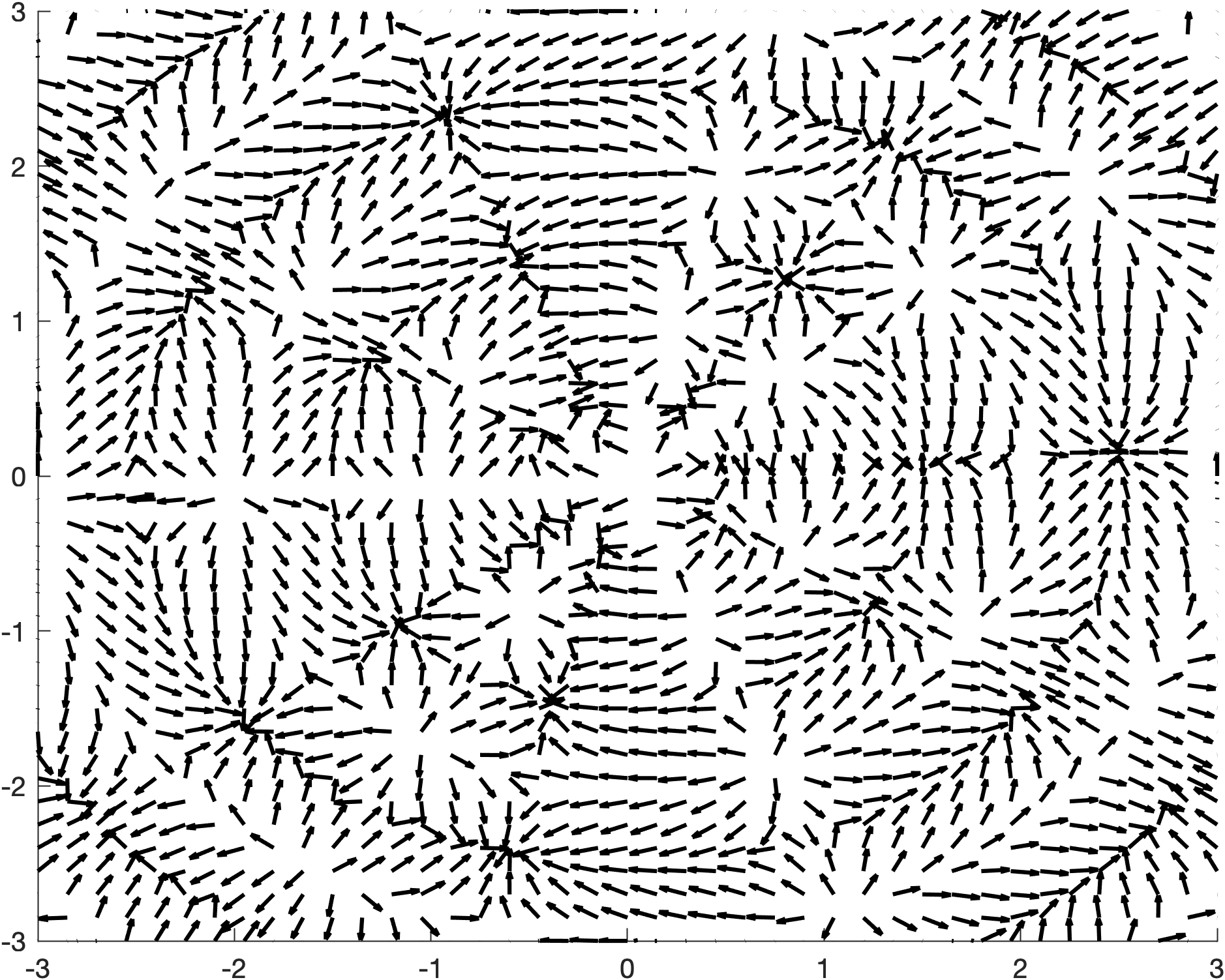}
        \caption{$t=0.66$}
        %\label{fig:ex_radial}
    \end{subfigure}
  \caption{(Normalized) Exemplary vector fields $X$ at different times\\
  $f(\theta)=2\cdot \cos(10 \pi \theta) \cdot \sin (4 \pi \theta) $  \\
    $g(r) = e^r\cdot \sin(2 \pi r)$}\label{fig:ex_vf}
\end{figure}  
We consider $1$-periodic orbits of the vector field $X$, i.e.~solutions $z:S^1=\R/\Z \to\R^2$ of $ \dot{z}(t) = X_t\big(z(t)\big)$.
Note that, since $X_t(0)=0$ for all $t\in S^1$, the constant map $S^1\ni t\mapsto 0\in \R^2$ is a 1-periodic orbit of $X$. From now on let us denote \begin{equation}
Z:=f^{-1}(1)\subset S^1 \quad\text{and}\quad R:=g^{-1}(0)\subset \R_{\geq 0}.
\end{equation} 
The set $Z$ may be empty, but $0\in R$ by definition of $g$.

In the next two lemmas, we classify all 1-periodic orbits of $X$.
\begin{lemma}
\label{lem:existence-of-periodic-orbits}
For  $\alpha \in Z=f^{-1}(1)$ and $\rho \in R=g^{-1}(0)$, the map
\begin{align}
    &z_{\alpha,\rho}: S^1 \to \R^2 \nonumber \\
    &z_{\alpha,\rho}(t)=\rho \cdot e^{2 \pi i (\alpha + t)}
\end{align}
is a 1-periodic orbit of the vector field \eqref{eq:vector-field}.
\end{lemma}

\begin{proof}
We include the straightforward calculation here for the convenience of the reader. For $\rho=0$ the map $z_{\alpha,\rho=0}=0$ is the constant orbit. Otherwise, 
\begin{align*}
    X_t\big( z_{\alpha,\rho}(t) \big)
    &= \underbrace{g(\rho)}_{=0} e^{2\pi i (\alpha+t)} + \underbrace{f(\alpha + t - t)}_{=1}  2\pi i  \rho  e^{2\pi i (\alpha+t)} \\
    &= 2\pi i  \rho  e^{2\pi i (\alpha+t)}\\
    &= \dot{z}_{\alpha,\rho}(t)
\end{align*}
shows the claim.
\end{proof}
 
\begin{lemma}
\label{lem:uniqueness-of-periodic-orbits}
Assume $Z \neq \emptyset$. Let $z:\R/k\Z \to \R^2$ be a $k$-periodic orbit for $k \in \N_{\geq 1}$.
Then there are $\alpha\in Z$ and $\rho \in R$ such that $z$ is a $k$-fold iteration of $z_{\alpha,\rho}$, that is
\begin{equation}\nonumber
\begin{aligned}
z(t)=z_{\alpha,\rho}(t\;\mathrm{mod}\;1).
\end{aligned}
\end{equation}
\end{lemma}

\begin{proof}
We, of course, assume that $z(t)\neq0$ for all $t$ since otherwise $z$ is the constant orbit $z\equiv0$. Thus, it is convenient to use polar coordinates $(r,\theta)\in\R_{\geq 0}\times S^1$, i.e.~$z=re^{2\pi i\theta}\in\C\setminus\{0\}$. In polar coordinates, the vector field $X$ is given by
\begin{equation}\nonumber
X_t(r, \theta)=\big(g(r), f(\theta-t)\big)
\end{equation}
and $\dot{z}(t) = X_t\big(z(t)\big)$ transforms into
\begin{align}\label{eq:vf_polar}
\left\{
\begin{aligned}
\dot r&= g(r)\\
\dot\theta&= f(\theta-t)
\end{aligned}
\right.
\end{align}
where the function $\theta(t)$ satisfies $\theta(t+k)=\theta(t)+\ell$ for all $t\in\R$ and some $\ell\in\Z$.

Since $r$ is $k$-periodic, the first equation in \eqref{eq:vf_polar} implies by uniqueness of solutions of ODEs that $r(t)\equiv\rho>0$ is constant and $g(\rho)=0$, that is $\rho\in R$. In particular, $z(t)=\rho e^{2\pi i\theta(t)}$. 

It remains to consider the second equation in \eqref{eq:vf_polar}, that is $\dot{\theta}(t)=f(\theta(t)-t)$ for all $t\in\R$. We recall that $t\mapsto \alpha +t$, $\alpha\in Z$, is a solution to the second equation in \eqref{eq:vf_polar}. 
Therefore, if $\theta(t_0)-t_0=\alpha\in Z$ for some $t_0$, then $\theta(t)=\alpha+t$ for all $t$ by uniqueness of solutions of ODEs.

Thus, we now assume that $\theta(t)-t\not\in Z$ for all $t$. In fact, we may assume that $0<\theta(t)-(\alpha+t)<1$ for all $t$ and all $\alpha\in Z$ since otherwise $e^{2\pi i\theta(t)}= e^{2\pi i(\alpha+t)}$ for some $t$ which again contradicts uniqueness. Now choose some $\alpha\in Z=\{f=1\}$ and consider
\begin{equation}\label{eqn:theta_alpha_t}
\begin{aligned}
\frac{\rd}{\rd t}(\theta(t)-(\alpha+t))&=\dot\theta(t)-1\\
&=f(\theta(t)-t)-f(\alpha+t-t)\\
&=f(\theta(t)-t)-1.
\end{aligned}
\end{equation}
By our assumption $\theta(t)-t\not\in Z$ we know that $f(\theta(t)-t)-1\neq0$. Let us assume for now that $f(\theta(t)-t)-1>0$. Since $\theta(t)$ satisfies $\theta(t+k)=\theta(t)+\ell$ and since $f$ is a periodic function there exists $\delta>0$ with $f(\theta(t)-t)-1\geq\delta>0$ for all $t\in\R$ and we conclude from \eqref{eqn:theta_alpha_t} that the difference $\theta(t)-(\alpha+t)$ grows at least as fast as $\delta t$, i.e.~cannot be bounded from above. If $f(\theta(t)-t)-1<0$ we similarly conclude that $\theta(t)-(\alpha+t)$ shrinks at least as fast as $-\delta t$ and cannot be bounded from below.

We conclude that our assumption $\theta(t)-t\not\in Z$ is wrong and therefore $\theta(t)=\alpha +t$ for some $\alpha\in Z$. In particular, $\ell=k$, that is $\theta(t+k)=\theta(t)+k$ and therefore, $z(t)$ is the $k$-fold cover of $z_{\alpha,\rho}$.
\end{proof}

To see when Theorem \ref{thm:mytheorem} applies to the explicit family $z_{\alpha,\rho}$ of 1-periodic orbits we first recall the definition of non-degeneracy.

\begin{definition}\label{def:non_degenerate}
    Let $x: S^1 \to \R^n$ be a $1$-periodic orbit of a (time-dependent) vector field $X$ on $\R^n$. Denote by $\Phi_t : \R^n \to \R^n$, $t\in\R$, the flow of $X$. The periodic orbit $x$ is called \textit{non-degenerate} if the linearized time-$1$-map $\rd \Phi_1\big( x_0(0)\big) :\R^n \to \R^n$ does not have $1$ as an eigenvalue.
\end{definition}

The following is an alternative characterization of non-degeneracy in terms of the vector field.

\begin{lemma}[{\cite[Lemma 6.6]{AlbersSeifert22}}]
\label{lem:non-deg-alternative}
Let $x : S^1 \rightarrow \mathbb{R}^n$ be a 1-periodic orbit of $X$. For $t\in S^1$ we consider the linear map
\begin{align*}
A(t) := -\rd X_t(x(t))^T : \mathbb{R}^n \longrightarrow \mathbb{R}^n.
\end{align*}
and denote by $Y: \R \to \R^{n\times n}$  the fundamental system for $A:S^1\to\R^{n\times n}$, that is, the solution of
\begin{align}
\begin{cases}  \frac{\rd}{\rd t} Y(t) = A(t) \cdot Y(t) \\[.5ex] Y(0) = \mathds{1}. \end{cases}
\label{eq:def-fundamental-system}
\end{align}
Then 
\begin{align*}
\rd\Phi^1_X(x(0))= \big(Y(1)^T\big)^{-1}
\end{align*}
holds and, in particular, $x$ is non-degenerate if and only if $Y(1)$ does not have $1$ as an eigenvalue.
\end{lemma}

Now we can characterize when an orbit $z_{\alpha,\rho}$ is non-degenerate.
 
 \begin{lemma}
 \label{lem:non-degeneracy}
  Fix $\alpha\in Z$ and $0\neq \rho\in R$. The $1$-periodic orbit $z_{\alpha,\rho}$ is non-degenerate if and only if $f'(\alpha)\neq 0$ and $g'(\rho) \neq 0$.
 \end{lemma}
 
 \begin{proof}
    We will use the alternative characterization of non-degeneracy from Lemma \ref{lem:non-deg-alternative}. Since the notion of non-degeneracy is coordinate independent we work again in polar coordinates. From $X_t(r, \theta)=\left(g(r), f(\theta-t)\right)$ together with $z_{\alpha,\rho}(t)=(\rho,\alpha+t)$ we conclude
\begin{align*}
A(t):=-\rd X_t\left(z_{\alpha,\rho}(t)\right)^T
        = \begin{pmatrix}
            -g'(\rho) & 0 \\
            0 & - f'(\alpha)
        \end{pmatrix}.
\end{align*}
Since $A(t)$ is time-independent (and diagonal) its fundamental system is simply given by $Y(t)=e^{tA}$, i.e.
\begin{align}
        Y(1)=e^A
        = \begin{pmatrix}
            e^{-g'(\rho)} & 0 \\
    0 & e^{-f'(\alpha)}
        \end{pmatrix}.
        \label{eq:Y}
\end{align}
Therefore, $1$ is not an eigenvalue of $Y(1)$ if and only  if $f'(\alpha)\neq 0$ and $g'(\rho) \neq 0$.
 \end{proof}

\section{1-periodic delay orbits}

Now we consider the delay differential equation
\begin{align}
    \label{eq:DDE}
    \dot{z}(t) = X_{t}\big(z(t-\tau)\big).
\end{align}
where $\tau\in\R$ and $X$ is as in \eqref{eq:vector-field}, that is,
\begin{align}
    &X: S^1 \times \R^2 \longrightarrow \R^2 \nonumber \\
    &X_t(z)=g(\Vert z\Vert) \cdot \frac{z}{\Vert z\Vert} + f\left(\frac{\arg(z)}{2 \pi} - t\right) \cdot 2\pi i   z\;.\nonumber
\end{align}
Since we require $g(0)=0$, the constant orbit $z\equiv 0$ also solves \eqref{eq:DDE} for every value of $\tau$. As opposed to the corresponding ordinary differential equation it is not clear (to us) how to construct or read off (periodic) delay orbits from the portrait of $X$, e.g.~Figure \ref{fig:ex_vf}. We point out that, again opposed to solutions to the ODE, different delay orbits may agree at various times without agreeing globally. In particular, a delay orbit may run into the constant orbit $z\equiv 0$ and then out again.

Let us fix $\alpha\in Z=\{f=1\}$ and $0\neq \rho\in R=\{g=0\}$. If $f'(\alpha) \neq 0$ and  $g'(\rho) \neq 0$, then the orbit $z_{\alpha,\rho}(t)=\rho e^{2 \pi i (\alpha + t)}$ defined in Lemma~\ref{lem:existence-of-periodic-orbits} is non-degenerate by Lemma~\ref{lem:non-degeneracy}. Thus, Theorem~\ref{thm:mytheorem} asserts that, for small enough $\tau\in\R$, there is a unique $1$-periodic solution $z_{\alpha,\rho,\tau}$ of the delay equation \eqref{eq:DDE} which is $\cC^\infty$-close to the given orbit $z_{\alpha,\rho}$. At this point, however, we do not know what $z_{\alpha,\rho,\tau}$ looks like. Moreover, if $f'(\alpha)=0$ or $g'(\rho)=0$ then $z_{\alpha,\rho}$ is a degenerate orbit and we cannot apply Theorem~\ref{thm:mytheorem}.

In Lemma \ref{lem:delay-orbits} below we will show that the loop
\begin{align*}
z_{\alpha,\rho,\tau}(t)=r_{\rho,\tau}\cdot e^{2 \pi i \cdot (t_{\alpha,\tau}+\tau+t)},
\end{align*}
where $r_{\rho,\tau}\in\R_{\geq 0}$ and $t_{\alpha,\tau}\in S^1$ are solutions to
\begin{align}
    f(t_{\alpha,\tau}) &= \cos(2 \pi \tau) \label{eq:equ-for-t-alpha-tau} \\
    g(r_{\rho,\tau}) &= -2\pi r_{\rho,\tau}  \sin(2 \pi \tau), \label{eq:equ-for-r-rho-tau}
\end{align}
indeed solves the delay equation \eqref{eq:DDE}. We point out that for $\tau=0$ the choice of $t_{\alpha,0}=\alpha\in Z$ and $r_{\rho,\tau}=\rho\in R$ recovers the (locally unique) solution $z_{\alpha,\rho}$ to the corresponding ODE. In accordance with Theorem \ref{thm:mytheorem} we expect to find delay orbits $\cC^\infty$-close to $z_{\alpha,\rho}$, at least for small delay $\tau$. This translates into $t_{\alpha,0}$ being close to $\alpha$ and $r_{\rho,\tau}$ being close to $\rho$.

\begin{lemma}
\label{lem:delay-orbits}
    Let $t_{\alpha,\tau}\in S^1$ and $r_{\rho,\tau}\in\R_{\geq 0}$ be solutions of \eqref{eq:equ-for-t-alpha-tau} and \eqref{eq:equ-for-r-rho-tau}, respectively. Then the loop
    \begin{align}
        & z_{\alpha,\rho,\tau}:S^1 \rightarrow \C \nonumber \\
        & z_{\alpha,\rho,\tau}(t)= r_{\rho,\tau} \cdot e^{2 \pi i \cdot (t_{\alpha,\tau}+\tau+t)}
        \label{eq:def-of-delay-orbit}
    \end{align}
    is a $1$-periodic $\tau$-delay orbit, i.e.~it solves the delay equation \eqref{eq:DDE} with delay $\tau$.
\end{lemma}

\begin{remark}
In Lemma \ref{lem:delay-orbits} we do not require the delay $\tau$ to be small. We point out, however, that, depending on $f,g,\rho,\alpha$ and $\tau$, solutions $t_{\alpha,\tau}\in S^1$ of \eqref{eq:equ-for-t-alpha-tau} resp.~$r_{\rho,\tau}\in\R_{\geq 0}$ of \eqref{eq:equ-for-r-rho-tau} may neither exist nor be unique.
Moreover, as long as $r_{\rho,\tau}$ and $t_{\alpha,\tau}$ are solutions of \eqref{eq:equ-for-t-alpha-tau} and \eqref{eq:equ-for-r-rho-tau} the loop $z_{\alpha,\rho,\tau}(t)= r_{\rho,\tau}  e^{2 \pi i \cdot (t_{\alpha,\tau}+\tau+t)}$ is a $\tau$-delay orbit. This does not require $r_{\rho,\tau}$ being close to some $\rho\in R$ or $t_{\alpha,\tau}$ being close to some $\alpha \in Z$. In fact, $Z=\emptyset$ or $R=\{0\}$ is allowed in Lemma \ref{lem:delay-orbits}.
\end{remark}

\begin{proof}[Proof of Lemma \ref{lem:delay-orbits}]
This is again a direct verification. Let us recall that
\begin{align*}\nonumber
X_t(z)=g(\Vert z\Vert)  \frac{z}{\Vert z\Vert} + f\left(\frac{\arg(z)}{2 \pi} - t\right)  2\pi i   z
\end{align*}
and
\begin{align*}\nonumber
z_{\alpha,\rho,\tau}(t)= r_{\rho,\tau}  e^{2 \pi i  (t_{\alpha,\tau}+\tau+t)}.
\end{align*}
We conclude immediately that
\begin{align*}\nonumber
\dot{z}_{\alpha,\rho,\tau}= 2 \pi i  z_{\alpha,\rho,\tau}.
\end{align*}
Using \eqref{eq:equ-for-t-alpha-tau} and \eqref{eq:equ-for-r-rho-tau}, we find that
    \begin{align*}
        X_t\big(z_{\alpha,\rho,\tau}(t-\tau)\big)
        &= g(r_{\rho,\tau})  e^{2 \pi i \big(t_{\alpha,\tau}  + \tau + (t - \tau)\big)} \\
        &\phantom{=}\quad + f\Big(\big(t_{\alpha,\tau} + \tau +(t  -\tau)\big) - t\Big)  2 \pi i r_{\rho,\tau}  e^{2 \pi i\big(t_{\alpha,\tau}  + \tau + (t - \tau)\big)}\\
        &= \underbrace{g(r_{\rho,\tau})}_{-2\pi r_{\rho,\tau}  \sin(2 \pi \tau)} e^{2 \pi i (t_{\alpha,\tau} +t)} + \underbrace{f(t_{\alpha,\tau})}_{ \cos(2 \pi \tau)} 2 \pi i r_{\rho,\tau}  e^{2 \pi i(t_{\alpha,\tau} +t)}\\
        &= 2\pi ir_{\rho,\tau}\big(\cos(2\pi \tau)+i \sin(2 \pi \tau)\big)e^{2 \pi i(t_{\alpha,\tau} +t)}\\
        &= 2\pi ir_{\rho,\tau}e^{2\pi i\tau}e^{2 \pi i(t_{\alpha,\tau} +t)}\\
        &= 2\pi ir_{\rho,\tau}e^{2 \pi i(t_{\alpha,\tau} +\tau +t)}\\
        &= 2 \pi i  z_{\alpha,\rho,\tau}(t)\\
        &=\dot{z}_{\alpha,\rho,\tau}(t),
    \end{align*}
showing that $z_{\alpha,\rho,\tau}$ is indeed a $1$-periodic $\tau$-delay orbit.
\end{proof}

\begin{remark}
As opposed to Lemma \ref{lem:uniqueness-of-periodic-orbits} we do not expect that there is a similar uniqueness statement for periodic delay orbits of the vector field $X$. 
\end{remark}

\section[Analyzing the families]{Analyzing the family $z_{\alpha,\rho,\tau}$}

We will now analyze the delay orbits $z_{\alpha,\rho,\tau}(t)= r_{\rho,\tau}  e^{2 \pi i  (t_{\alpha,\tau}+\tau+t)}$ arising for different choices of $\alpha\in Z=\{f=1\}$ and $0\neq\rho\in R=\{g=0\}$ and corresponding solutions $t_{\alpha,\tau}$ of \eqref{eq:equ-for-t-alpha-tau} and  $r_{\rho,\tau}$ of \eqref{eq:equ-for-r-rho-tau}.

We note that the 1-periodic delay orbit $z_{\alpha,\rho,\tau}$ moves with constant angular speed on the circle of radius $r_{\rho,\tau}$. In particular, we can recover $z_{\alpha,\rho,\tau}$ from a single value, e.g.~$z_{\alpha,\rho,\tau}(0)$. Therefore, we choose to represent a family $(z_{\alpha,\rho,\tau})_{\tau\in[\tau_-,\tau_+]}$ of 1-periodic delay orbits by the curve 
\begin{align}
    \label{eq:def-of-curve}
   \gamma_{\alpha,\rho}:  [\tau_-,\tau_+] &\longrightarrow \R^2 \nonumber \\
    \tau &\longmapsto z_{\alpha,\rho,\tau}(0)=r_{\rho,\tau}  e^{2 \pi i  (t_{\alpha,\tau}+\tau)}.
\end{align}
From this explicit expression it is clear that $\gamma_{\alpha,\rho}$ is smooth if $t_{\alpha,\tau}$ and $r_{\rho,\tau}$ smoothly depend on $\tau$. 

\subsection{The non-degenerate case} \label{subsec:non-deg-case}

We recall from Lemma \ref{lem:non-degeneracy} that the 1-periodic orbit $z_{\alpha,\rho}$ of $X$ is non-degenerate if and only if $f'(\alpha)\neq 0$ for $\alpha\in Z=\{f=1\}\subset S^1$ and $g'(\rho)\neq 0$ for $0\neq \rho\in R=\{g=0\}\subset\R_{\geq0}$. In fact, non-degeneracy of $z_{\alpha,\rho}$, i.e.~$f'(\alpha)\neq 0$ and $g'(\rho)\neq 0$, allows us to locally uniquely solve equations \eqref{eq:equ-for-t-alpha-tau} and  \eqref{eq:equ-for-r-rho-tau} as follows. For that first choose a local inverse 
\begin{align*}
f^{-1}_{\alpha}:I_{\alpha}\to N_{\alpha}
\end{align*} 
from an interval $I_{\alpha}\subset \R$ containing $1$ onto a neighborhood $N_{\alpha}\subset S^1 $ of $\alpha$. In particular, $f_\alpha^{-1}(1)=\alpha$. Then, for $\tau$ small enough (more precisely, whenever $\cos(2\pi\tau)\in I_{\alpha}$), we set
\begin{align} \label{eq:def-of-t-alpha-tau}
    t_{\alpha,\tau}:=f^{-1}_{\alpha}\big(\cos(2\pi\tau)\big).
\end{align}
By construction, $t_{\alpha,\tau}$ solves \eqref{eq:equ-for-t-alpha-tau}. Next, we abbreviate
\begin{align*}
\tilde g(r):=\frac{g(r)}{r}:\R_{>0}\to\R
\end{align*}
and observe that the assumption $g'(\rho)\neq0$ for $\rho>0$ implies
\begin{align*}
\tilde{g}'(\rho) = \frac{g'(\rho)\rho-g(\rho)}{\rho^2} > 0.
\end{align*}
Again, we choose a local inverse 
\begin{align*}
\tilde{g}^{-1}_{\rho}: I_{\rho}\to N_{\rho}
\end{align*} 
from an interval $I_{\rho}\subset \R$ containing $0$ to a  neighborhood $N_{\rho}\subset\R_{>0}$ of $\rho$, in particular $\tilde{g}^{-1}_{\rho}(0)=\rho$, and set 
\begin{align} \label{eq:def-of-r-rho-tau}
    r_{\rho,\tau}:=\tilde{g}^{-1}_{\rho} \big(-2\pi \sin(2\pi\tau)\big)
\end{align}
which, by construction, $r_{\rho,\tau}$ solves \eqref{eq:equ-for-r-rho-tau}.

For convenience let us assume that $1\in I_{\alpha}\subset\R$ resp.~$0\in I_{\rho}\subset\R$ are the maximal intervals on which the inverses $f^{-1}_{\alpha}$ resp.~$\tilde{g}^{-1}_{\rho}$ are well-defined.

\begin{figure}[h]
\captionsetup[subfigure]{}%{justification=centering}
  \begin{subfigure}[t]{.47\textwidth}
    \centering
    \includegraphics[width=\linewidth]{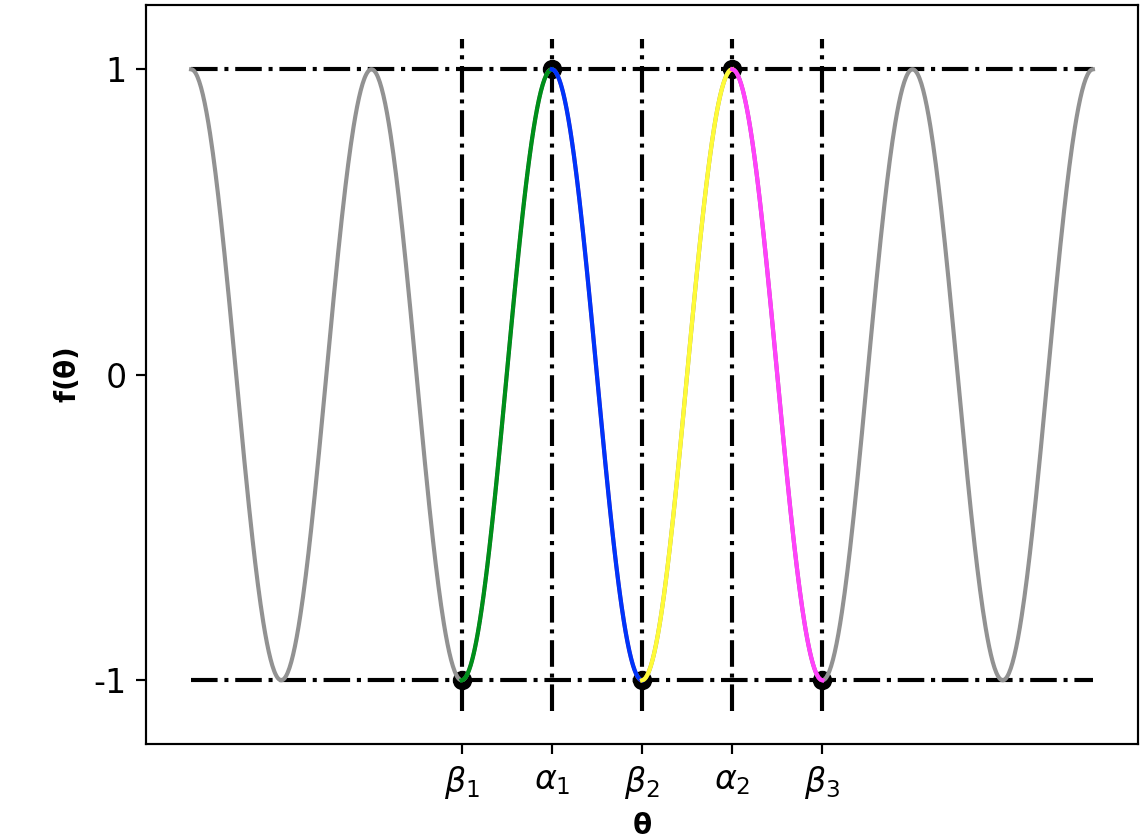}
    \caption{If $f$ takes a maximum in $\alpha$ (or a minimum in $\beta$), there is a choice of local inverse.}
    \label{fig:deg_min}
  \end{subfigure}
  \hfill
  \begin{subfigure}[t]{.47\textwidth}
    \centering
    \includegraphics[width=\linewidth]{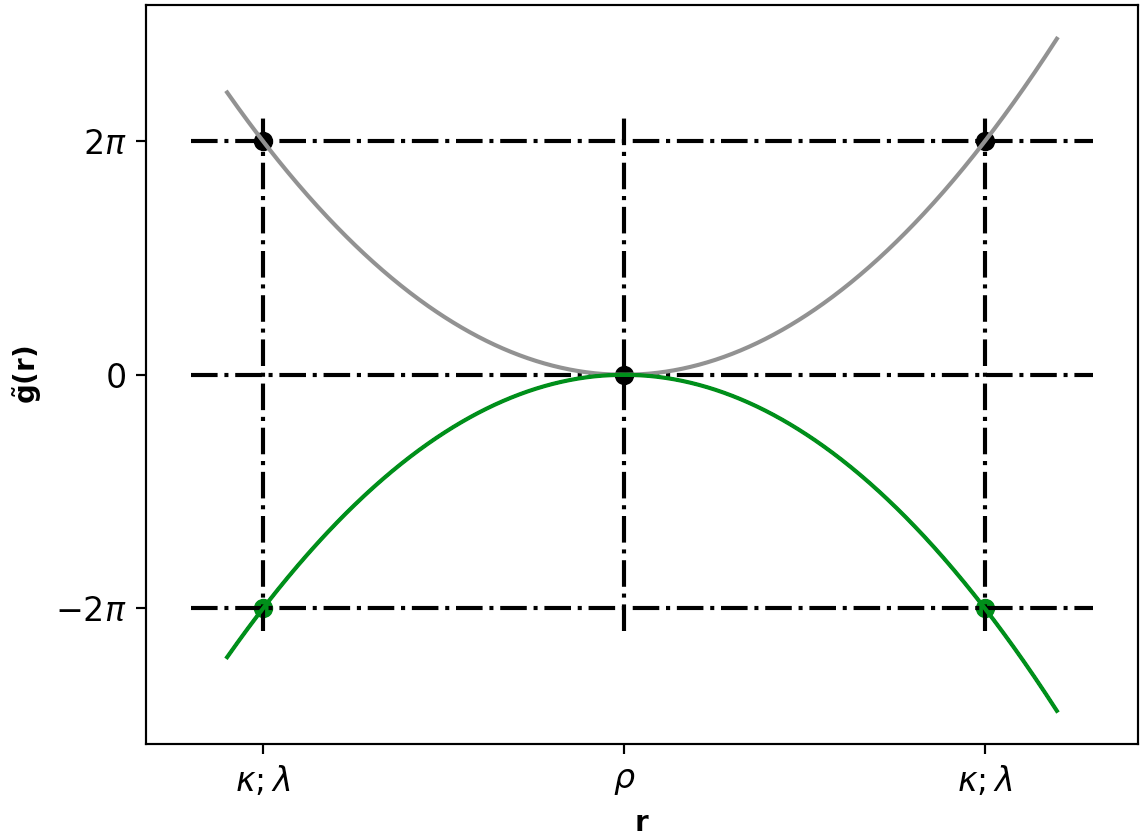}
    \caption{If $\tilde{g}$ takes a minimum or maximum in $\rho$, we can define $r_{\rho,\tau}$ only for $\tau<0$ or only  for $\tau>0$, respectively.}
    \label{fig:deg_g}
  \end{subfigure}
  \caption{The degenerate case.}\label{fig:deg_cases}
\end{figure}

\begin{remark} \label{rem:collection_of_results_in_nondeg_case}
We collect a few conclusions.
\begin{enumerate}
\item If $\tau$ (not necessarily small) is such that $\cos(2\pi\tau)\in I_{\alpha}$ and $-2\pi\sin(2\pi\tau)\in I_{\rho}$, we define $t_{\alpha,\tau}$ resp.~$r_{\rho,\tau}$ by \eqref{eq:def-of-t-alpha-tau} resp.~\eqref{eq:def-of-r-rho-tau}. In particular, we then obtain a 1-periodic $\tau$-delay orbit $z_{\alpha,\rho,\tau}(t)= r_{\rho,\tau}  e^{2 \pi i  (t_{\alpha,\tau}+\tau+t)}$.

Moreover, the so defined maps $\tau\mapsto t_{\alpha,\tau}$ and $\tau\mapsto r_{\rho,\tau}$ are continuous. They are smooth except for $\tau$ with $f'(t_{\alpha,\tau})=0$ or $\tilde{g}'(r_{\rho,\tau})=0$.

\item Sometimes $[-1,1]\subset I_{\alpha}\subset\R$ and $[-2\pi,2\pi]\subset I_{\rho}\subset\R$, see Figure (\ref{fig:con_delay_families}) and Figure (\ref{fig:ex_cusps_vf}) below. In this case the maps $\tau\mapsto t_{\alpha,\tau}$ and $\tau\mapsto r_{\rho,\tau}$ are defined for all $\tau\in\R$. Moreover, they are $1$-periodic in $\tau$. Hence, in this case the $\tau$-delay orbits $z_{\alpha,\rho,\tau}$, $\tau\in\R$, actually form an $S^1$-family. 

More precisely, as $\tau$ runs through $\R$ the value $t_{\alpha,\tau}$ oscillates between $\alpha=f_\alpha^{-1}(1)$, for $\tau\in\Z$, and $\beta:=f_\alpha^{-1}(-1)$, for $\tau\in\tfrac12+\Z$, and similarly the value $r_{\rho,\tau}$ oscillates between $\tilde{g}^{-1}_{\rho}(-2\pi)$, for $\tau\in\frac14+\Z$, and $\tilde{g}^{-1}_\rho(2\pi)$, for $\tau\in\frac34+\Z$.

Let us discuss possible cases a bit further.

\begin{enumerate}  
\item We have $[-1,1]\subset I_{\alpha}$ if and only if there is some $\beta\in S^1$ with $f(\beta)=-1$ and such that $f$ is strictly monotone on a segment in $S^1$ from $\alpha$ to $\beta$. In this case, the set $N_\alpha$ from above is exactly this segment from $\alpha$ to $\beta$. 

\item If we, in addition, know that $f'\big|_{N_\alpha}\neq0$, then the local inverse $f_\alpha^{-1}\big|_{[-1,1]}:[-1,1]\stackrel{\cong}{\to}N_\alpha$ is smooth, implying that $\tau\mapsto t_{\alpha,\tau}$ is smooth as well.

\item On the other hand, if $\beta$ is a non-degenerate local minimum (i.e.\ $f'(\beta)=0$ and $f''(\beta)>0$) with $f(\beta)=-1$ then we still have $[-1,1]\subset I_{\alpha}$ and we could proceed as above. On the other hand, we could also choose to continue the map $\tau\mapsto t_{\alpha,\tau}$ by switching from the local inverse $f_\alpha^{-1}$ at $\beta=f_\alpha^{-1}(-1)$ to the local inverse of $f$  on ``the other side'' of $\beta$, see Figure (\ref{fig:deg_min}). 

\item Similarly, we have $[-2\pi,2\pi]\subset I_{\rho}$ if and only if there are $\kappa, \lambda\in \R_{> 0}$ with $\tilde{g}(\kappa)=2\pi$, $\tilde{g}(\lambda)=-2\pi$ such that either $\kappa<\rho<\lambda$ or $\lambda<\rho<\kappa$ and $\tilde{g}$ is strictly monotone between $\kappa$ and $\lambda$. The set $N_{\rho}$ is then exactly the interval between $\kappa$ and $\lambda$. Again, $\tilde{g}'\big|_{N_\rho}\neq0$ implies that $\tilde{g}^{-1}_{\rho}\big|_{[-2\pi,2\pi]}:[-2\pi,2\pi]\stackrel{\cong}{\to} N_{\rho}$ is smooth and thus so is $\tau\mapsto r_{\rho,\tau}$.
A similar discussion as above for $f_\alpha$ is possible.

\item We call an $S^1$-family $z_{\alpha,\rho,\tau}$ \textit{smooth} if both local inverse $f_\alpha$ and $\tilde{g}_\rho$ are smooth functions, i.e.~if $f'\big|_{N_\alpha}\neq0$ and $\tilde{g}'\big|_{N_\rho}\neq0$ holds. This implies, that its associated curve $\gamma_{\alpha, \rho}$ is smooth.

\end{enumerate}

\item If $\tilde{g}(0):=g'(0)=\lim_{r\to 0}\frac{g(r)}{r}$ exists, then we can allow for $0\in N_{\rho}$, and hence $r_{\rho,\tau}$ may become $0$, namely if $\tau_*\in\frac14+\Z$ or $\tau_*\in\frac34+\Z$. In this case the family of $\tau$-delay orbits $z_{\alpha,\rho,\tau}$ contains the constant delay orbit $z\equiv 0$ at $\tau=\tau_*$. 

Without further assumptions on $g$, it is not clear that $\tilde{g}$ extends continuously to $0$ and therefore it is also not clear how the local inverse $\tilde{g}_\rho^{-1}$ behaves as the values $r_{\rho,\tau}$ approach $0$. In particular, we can only continue the family $z_{\alpha,\rho,\tau}$ beyond $\tau=\tau_*$ if $\tau_*\in\frac14+\Z$ and $\tilde{g}(0)=-2\pi$ or if  $\tau_*\in\frac34+\Z$ and $\tilde{g}(0)=2\pi$.

We do point out, however, that, if $g(r)$ is $\cC^l$ in $r=0$, then $\tilde{g}(r)=\tfrac{1}{r}g(r)$ admits a $\cC^{l-1}$ extension to $r=0$, since $\displaystyle g(r)=r \int_0^1g'(sr)\rd s$.
\end{enumerate}
\end{remark}
\begin{figure}[h]
\captionsetup[subfigure]{justification=centering}
 \begin{subfigure}[t]{.49\textwidth}
   \centering
   \includegraphics[width=\linewidth]{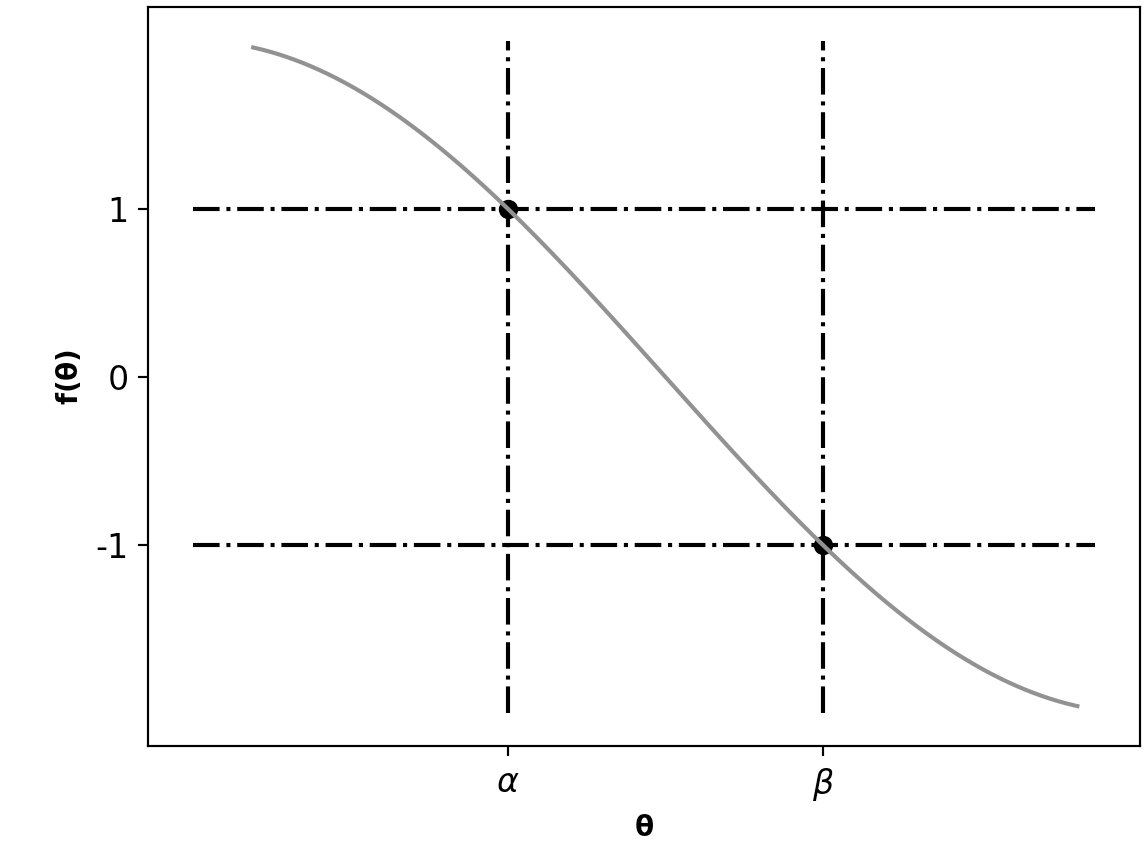}
   \caption{Construction of $t_{\alpha,\tau}$ for monotone $f$}
   \label{fig:con_f}
 \end{subfigure}
 \hfill
 \begin{subfigure}[t]{.49\textwidth}
   \centering
   \includegraphics[width=\linewidth]{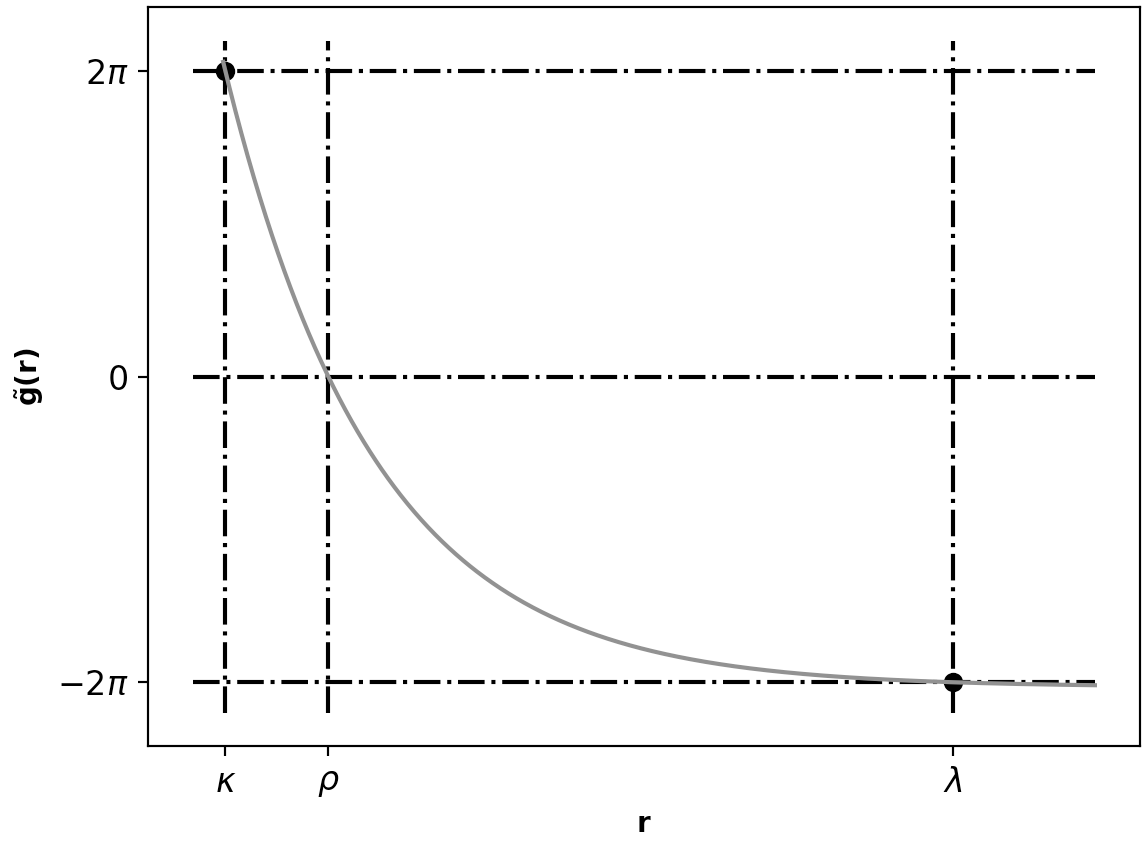}
   \caption{Construction of $r_{\rho,\tau}$ for monotone $\g$}
   \label{fig:con_g}
 \end{subfigure}
 \caption{Construction of $t_{\alpha,\tau}$ and $r_{\rho,\tau}$ for monotone $f$ and $\tilde{g}$}
 \label{fig:con_delay_families}
\end{figure}

\subsection{Cusps}

\begin{figure}[h]
\captionsetup[subfigure]{justification=centering}
    \begin{subfigure}[t]{0.4\linewidth}
        \centering
        \includegraphics[width=\linewidth,keepaspectratio]{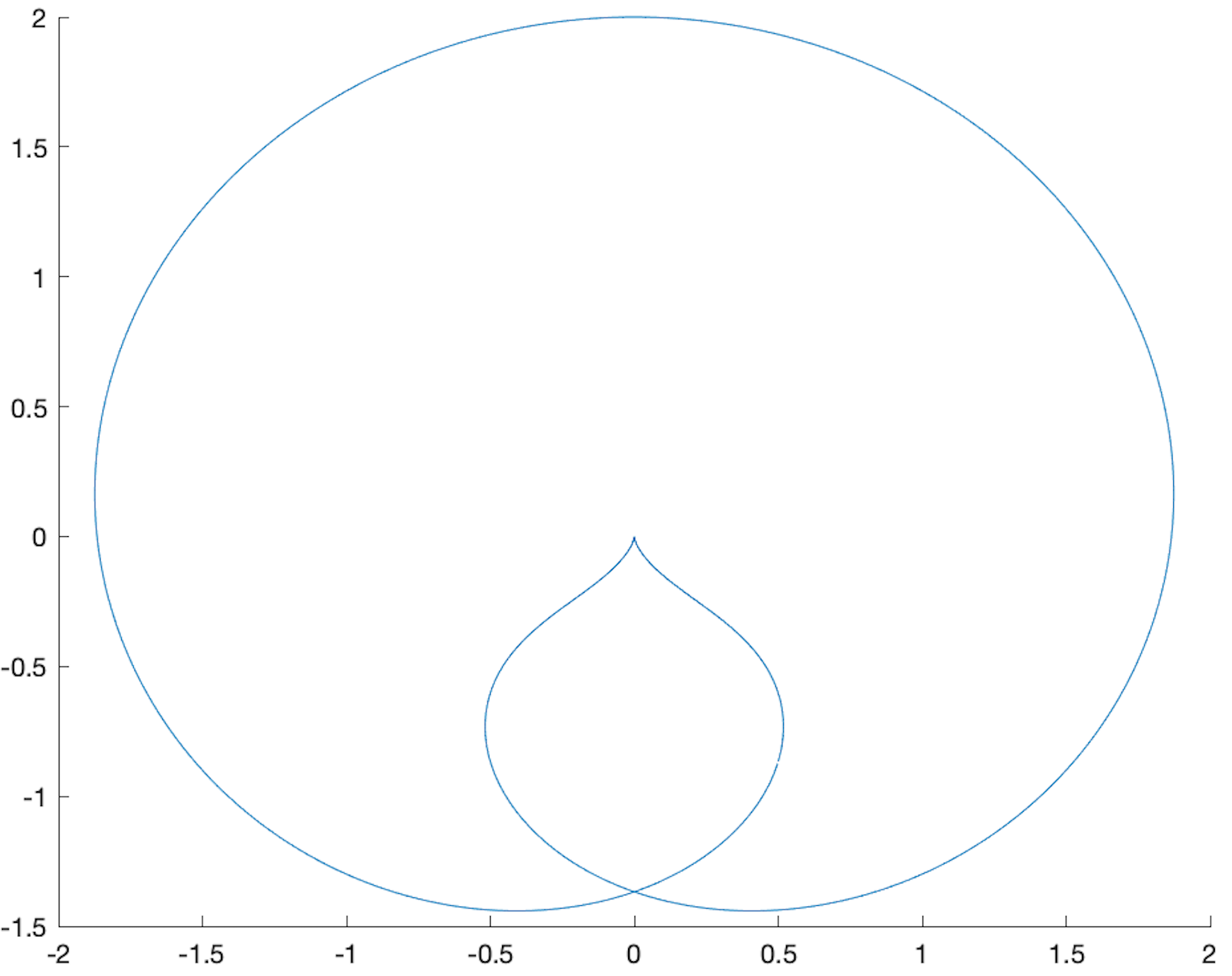}
         \caption{$f(\theta) = 3 \theta - 1.5$\\
         $g(r)= 2\pi x(x-1)$}
    \end{subfigure}
    \hfill
    \begin{subfigure}[t]{0.4\linewidth}
        \centering
         \includegraphics[width= \linewidth,keepaspectratio]{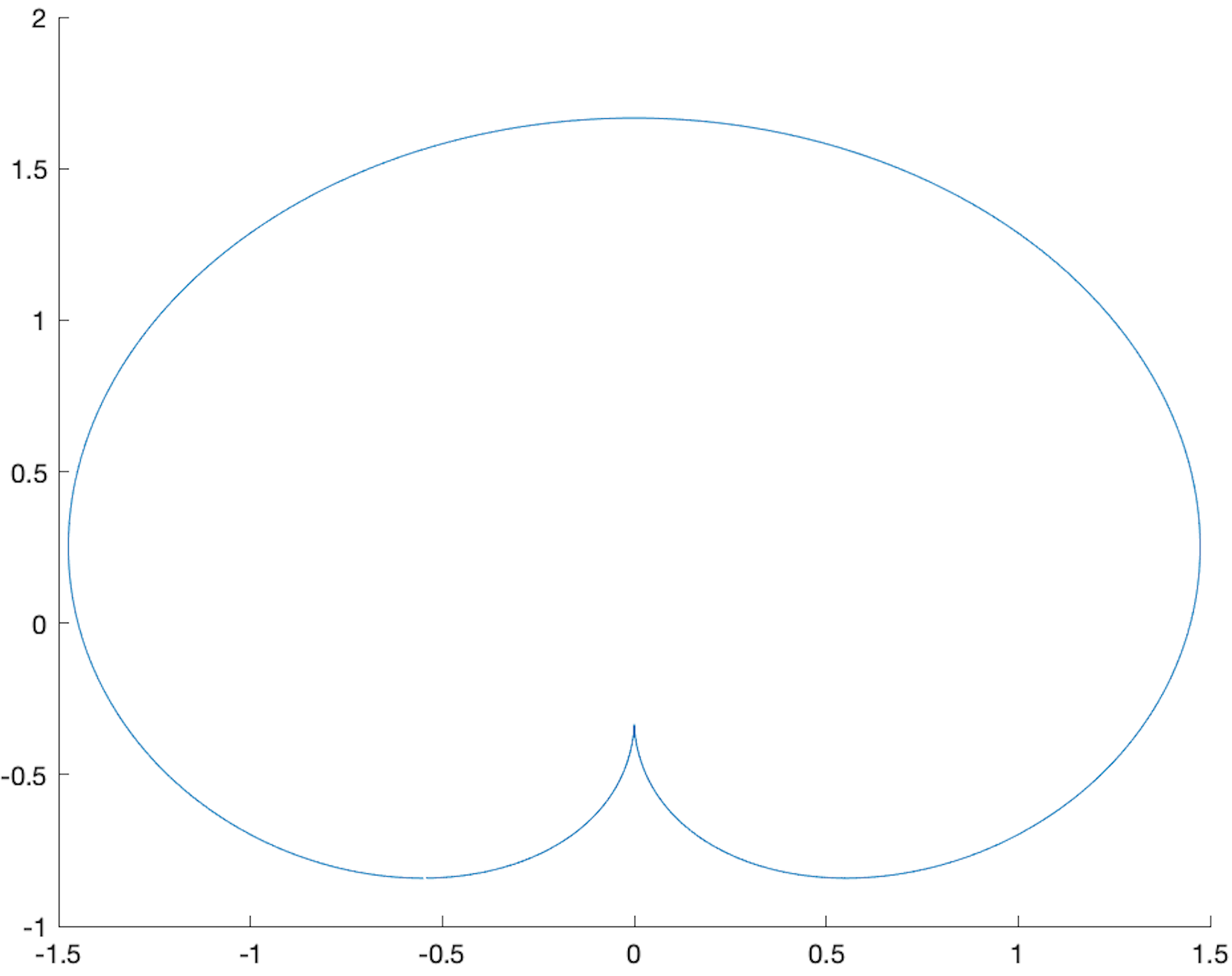}
         \caption{$f(\theta)= 2 \pi \theta - \pi$\\
         $g(r)= 3\pi x(x-1)$}
    \end{subfigure}
    
    \caption{On the left we have an example for a cusp due to condition (i) and on the right for condition (ii) of Lemma 4.2. }
    \label{fig:ex_cusps_one}
\end{figure}

\begin{figure}[h]
\captionsetup[subfigure]{justification=centering}
    \begin{subfigure}[t]{0.49\linewidth}
        \centering
        \includegraphics[width=\linewidth,keepaspectratio]{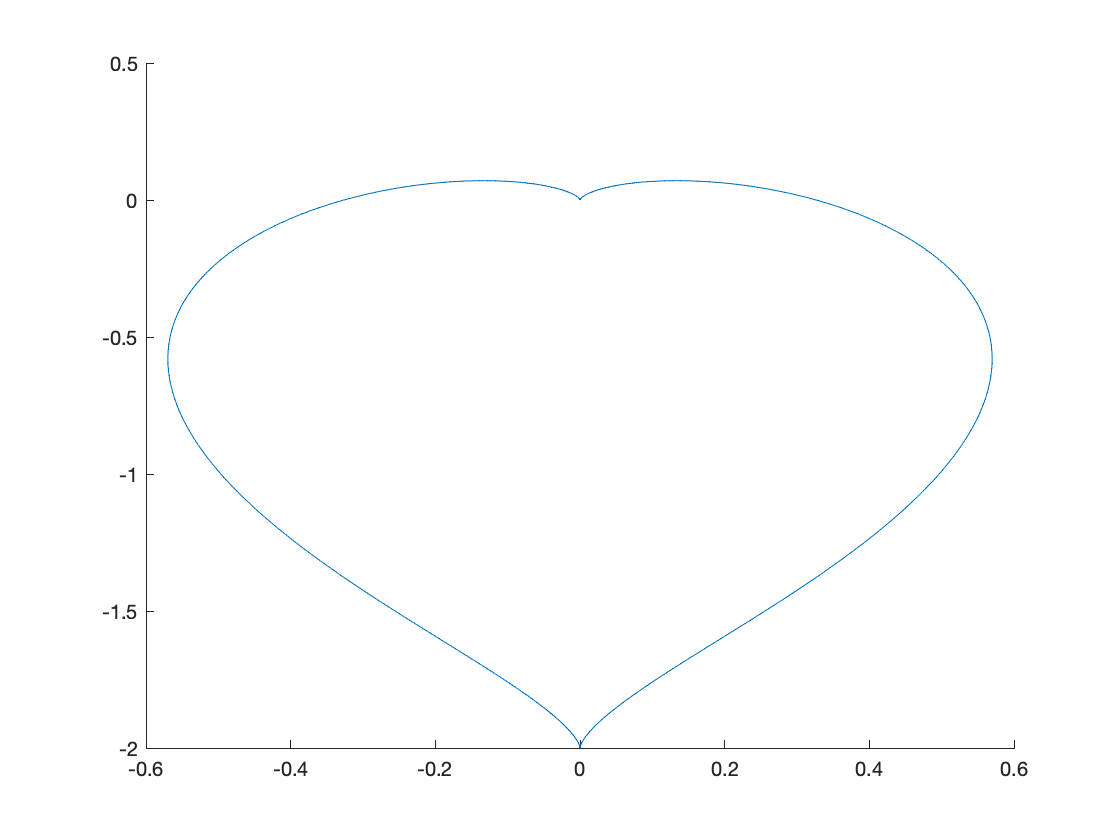}
         \caption{$f(\theta) = 2\pi \theta - \pi$\\
         $g(r)= -2\pi x(x-1)$}
         \label{fig:ex_df_cusps}
    \end{subfigure}
    \hfill
    \begin{subfigure}[t]{0.49\linewidth}
        \centering
         \includegraphics[width= \linewidth,keepaspectratio]{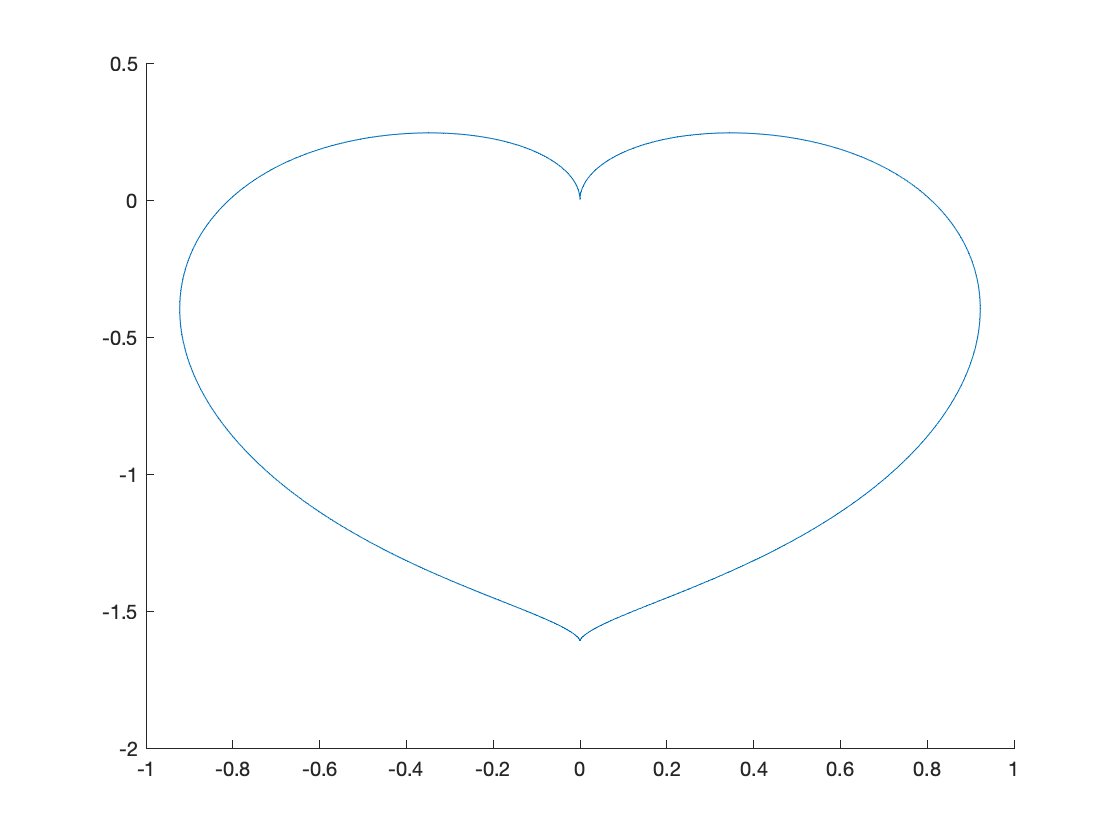}
         \caption{$f(\theta)= 2 \pi \theta - \pi$\\
         $g(r)= \frac{10}{5} (r^3 + 1)\sin(\frac{4}{5}\pi r)$}
         \label{fig:ex_df_cusps_2}
    \end{subfigure}   
    \caption{Examples of delay families with two cusps}
    \label{fig:ex_two_cusps}
\end{figure}

Recall that we represent families of delay orbits $(z_{\alpha,\rho,\tau}=r_{\rho,\tau}  e^{2 \pi i (t_{\alpha,\tau}+\tau)})_{\tau\in[\tau_-,\tau_+]}$ by the curve $\gamma_{\alpha,\rho}:  [\tau_-,\tau_+] \longrightarrow \R^2$ in the plane given by $\gamma_{\alpha,\rho}(\tau)=z_{\alpha,\rho,\tau}(0)$, see \eqref{eq:def-of-curve} for details. In example plots of the curve $\gamma_{\alpha,\rho}$ often exhibit isolated cusps-type singularities, see Figure  (\ref{fig:ex_cusps_one}) and (\ref{fig:ex_two_cusps}), meaning that if the derivative of $\gamma_{\alpha,\rho}$ has a zero at some $\tau_*$ then this is an isolated zero and there necessarily is a switch of direction of the normalized tangent vector to $\gamma_{\alpha,\rho}$ at $\tau_*$. We are going to explore this further now. For that, we assume in this section that the maps $\tau\mapsto t_{\alpha,\tau}$ and $\tau\mapsto r_{\rho,\tau}$ both are smooth.

\begin{lemma} \label{lem:app_cusp}
The curve $\gamma_{\alpha,\rho}$ has only isolated singularities. More precisely, the condition $\gamma_{\alpha,\rho}'(\tau_*)=0$ implies that $\tau_* \in \frac14+\Z$ or $\tau_*\in \frac34+\Z$. 

Conversely, assume $\tau_* \in \frac14+\Z$ or $\tau_*\in \frac34+\Z$. Then $\gamma_{\alpha,\rho}'(\tau_*)=0$ if one of the following holds:
\begin{enumerate}\renewcommand{\labelenumi}{\rm(\roman{enumi})}
\item $r_{\rho,\tau_*}= 0$,
\item $f'(t_{\alpha,\tau_*})= 2\pi$ in case $\tau_* \in \frac14+\Z$ resp.~$f'(t_{\alpha,\tau_*})= -2\pi$ in case $\tau_* \in \frac34+\Z$.
\end{enumerate}
\end{lemma}

\begin{figure}[h]
    \captionsetup[subfigure]{justification=centering}
    \begin{subfigure}[t]{0.49\linewidth}
        \centering
         \includegraphics[width=\linewidth,keepaspectratio]{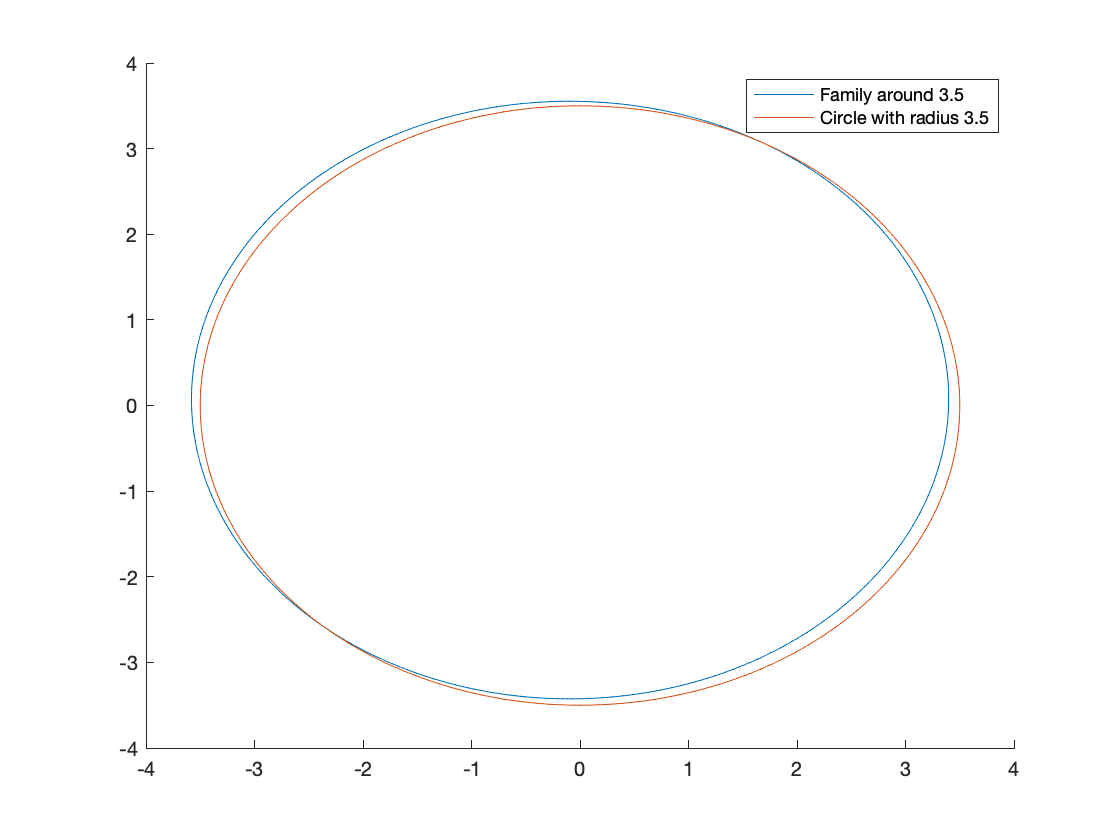}
         \caption{The image of the curve $\gamma_{\alpha_1,\rho_1}$ for $\alpha_1 = 0.1702 $ and $\rho_1=3.5$ in blue and for comparison a red circle of radius $3.5$.}
         \label{fig:ex_df_1}
    \end{subfigure}
    \hfill
    \begin{subfigure}[t]{0.49\linewidth}
        \centering
        \includegraphics[width=\linewidth,keepaspectratio]{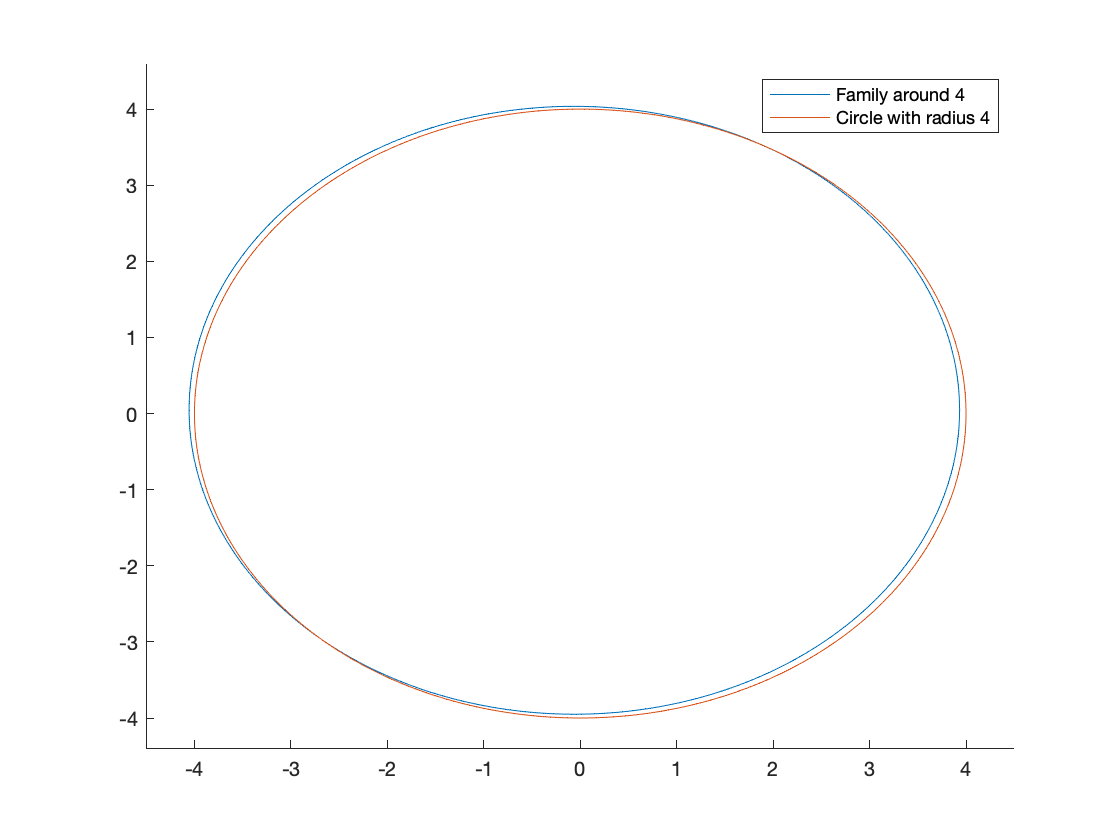}
         \caption{The image of the curve $\gamma_{\alpha_2,\rho_2}$ for $\alpha_2 = 0.6332 $ and $\rho_2=4$ in blue and for comparison a red circle of radius $4$.}
         \label{fig:ex_df_2}
    \end{subfigure}\\
%    \medskip
    \begin{subfigure}[t]{0.49\linewidth}
        \centering
        \includegraphics[width=\linewidth]{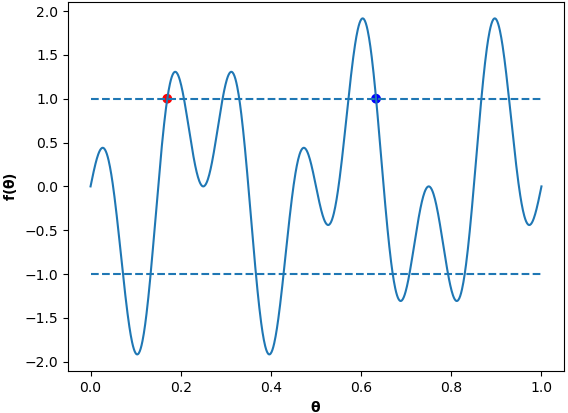}
        \caption{Plot of $f(\theta)=2 \cos(10 \pi \theta)  \sin (4 \pi \theta)$.\\ The dots mark the choices for $\alpha$ \\ (red: $\alpha_1$, blue: $\alpha_2$). }
        \label{fig:ex_angular}
    \end{subfigure}
    \hfill
    \begin{subfigure}[t]{0.49\linewidth}
        \centering
        \includegraphics[width=\linewidth]{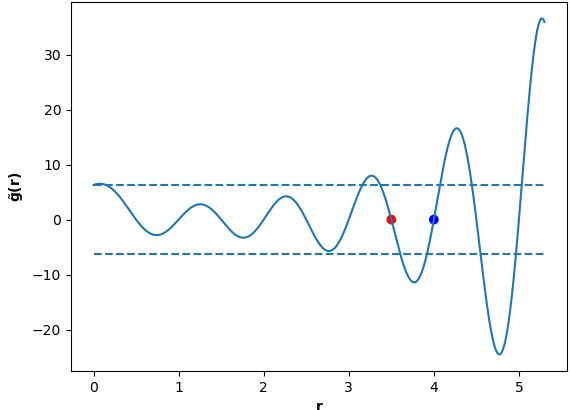}
        \caption{Plot of $\g(r)= \frac{1}{r} e^r \sin(2 \pi r)$. \\ The dots mark the choices for $\rho$ \\ (red: $\rho_1$, blue: $\rho_2$). }
        %\label{fig:ex_radial}
    \end{subfigure}
    \caption{Families of delay orbits for the example in Figure \ref{fig:example-vfield-a}. The figures show the smooth $S^1$ families we constructed in the discussion above. Moreover, the role and explicit choice of $\alpha$ and $\rho$ is depicted to show their influence on the delay families.}
    \label{fig:ex_cusps_vf}
\end{figure}

\begin{proof}
We first discuss the condition $\gamma_{\alpha,\rho}'(\tau_*) = 0$. From 
\begin{align}\label{eq:formula-derivative-of-curve}
\gamma_{\alpha,\rho}'(\tau_*) &= \frac{\rd}{\rd\tau}\Big|_{\tau=\tau_*}z_{\alpha,\rho,\tau}(0)\nonumber\\
&=  \frac{\rd}{\rd\tau} \Big|_{\tau=\tau_*}  r_{\rho,\tau}  e^{2 \pi i (t_{\alpha,\tau}+\tau)}  \nonumber\\
&=\left[ \left( \frac{\rd}{\rd\tau} \Big|_{\tau=\tau_*} r_{\rho,\tau} \right)+ r_{\rho,\tau_*}  2\pi i\left(\left(\frac{\rd}{\rd\tau} \Big|_{\tau=\tau_*} t_{\alpha,\tau}\right) +1 \right) \right] e^{2 \pi i (t_{\alpha,\tau_*}+\tau_*)}
\end{align}
we conclude that $\gamma_{\alpha,\rho}'(\tau_*) = 0$ is equivalent to
\begin{align}
        \label{eq:cusp-necessary-condition}
        \frac{\rd}{\rd\tau} \Big|_{\tau=\tau_*} r_{\rho,\tau}=0
        \qquad \text{ and } \qquad
        r_{\rho,\tau_*}  \left(\left(\frac{\rd}{\rd\tau} \Big|_{\tau=\tau_*} t_{\alpha,\tau}\right) +1\right)=0.
    \end{align}
Let us analyze these two equations individually. Recalling the definition \eqref{eq:def-of-r-rho-tau} of $r_{\rho,\tau}$ and that $\tilde g_\rho^{-1}$ is a local inverse for $\tilde g$ leads to
\begin{equation}\label{eqn:formula_d_dtau_r_rho_tau}
\frac{\rd}{\rd\tau} \Big|_{\tau=\tau_*} r_{\rho,\tau} = \frac{\rd}{\rd\tau} \Big|_{\tau=\tau_*} \tilde{g}_\rho^{-1}(-2\pi \sin(2 \pi \tau))= \frac{-4 \pi^2 \cos(2\pi \tau_*)}{\tilde{g}'(r_{\rho,\tau_*})}
\end{equation}
and we conclude that  $\frac{\rd}{\rd\tau} \big|_{\tau=\tau_*} r_{\rho,\tau}=0$ is equivalent to $\tau_* \in \frac14+\Z$ or $\tau_*\in \frac34+\Z$. In particular, we see that $\gamma_{\alpha,\rho}'(\tau_*)=0$ implies $\tau_* \in \frac14+\Z$ or $\tau_*\in \frac34+\Z$. 

The second equation in \eqref{eq:cusp-necessary-condition} holds if $r_{\rho,\tau_*}=0$ or if $\frac{\rd}{\rd\tau} \big|_{\tau=\tau_*} t_{\alpha,\tau}=-1$. As above the definition \eqref{eq:def-of-t-alpha-tau} of $t_{\alpha,\tau}$ together with the fact that $f_\alpha^{-1}$ is a local inverse to $f$ gives
\begin{equation}
\frac{\rd}{\rd\tau} \Big|_{\tau=\tau_*} t_{\alpha,\tau}= \frac{\rd}{\rd\tau} \Big|_{\tau=\tau_*} f_\alpha^{-1}(\cos(2 \pi \tau)) = \frac{- 2 \pi \sin(2 \pi \tau_*)}{f'(t_{\alpha,\tau_*})}.
\end{equation}
Now, if we assume that  $\tau_* \in \frac14+\Z$ resp.~$\tau_*\in \frac34+\Z$, then we see that $\frac{\rd}{\rd\tau} \big|_{\tau=\tau_*} t_{\alpha,\tau}=-1$ is equivalent to
\begin{equation}
 f'(t_{\alpha,\tau_*})= 2\pi \qquad \text{resp.} \qquad f'(t_{\alpha,\tau_*})= -2\pi
\end{equation}
proving the Lemma.
\end{proof}

\begin{proposition}
 Assume that $\gamma_{\alpha,\rho}'(\tau_*) = 0$. Then the curve $\gamma_{\alpha,\rho}$ necessarily has a cusp at $\tau_*$, that is
\begin{equation}\label{eq:cusp_condition}
\lim_{\tau\nearrow\tau_*}\frac{\gamma_{\alpha,\rho}'(\tau)}{\Vert\gamma_{\alpha,\rho}'(\tau)\Vert}=-\lim_{\tau\searrow\tau_*}\frac{\gamma_{\alpha,\rho}'(\tau)}{\Vert\gamma_{\alpha,\rho}'(\tau)\Vert}\neq0.
\end{equation}
I.e., the normalized tangent vector switches direction at $\tau_*$.
\end{proposition}

\begin{proof}
Since the maps $\tau\mapsto t_{\alpha,\tau}$ and $\tau\mapsto r_{\rho,\tau}$ are smooth, so is $z_{\alpha,\rho,\tau}=r_{\rho,\tau}  e^{2 \pi i (t_{\alpha,\tau}+\tau)}$ and we use first order Taylor approximation of $\gamma_{\alpha,\rho}'(\tau)=\frac{\rd}{\rd\tau}z_{\alpha,\rho,\tau}(0)$ near $\tau=\tau_*$, i.e.
\begin{equation}
\gamma_{\alpha,\rho}'(\tau_*+\varepsilon)= \gamma_{\alpha,\rho}'(\tau_*)+\gamma_{\alpha,\rho}''(\tau_*) \varepsilon+r(\varepsilon)\varepsilon=\gamma_{\alpha,\rho}''(\tau_*) \varepsilon+ r(\varepsilon)\varepsilon
\end{equation}
for some function $r(\varepsilon)$ with $\lim_{\varepsilon\to0}r(\varepsilon)=0$. Therefore, the desired condition \eqref{eq:cusp_condition} holds, if we can show that $\gamma_{\alpha,\rho}''(\tau_*)\neq0$. For this, we recall from the proof of Lemma \ref{lem:app_cusp} that
\begin{equation}
\gamma_{\alpha,\rho}'(\tau_*) =\left[ \left( \frac{\rd}{\rd\tau} \Big|_{\tau=\tau_*} r_{\rho,\tau} \right)+ r_{\rho,\tau_*}  2\pi i\left(\left(\frac{\rd}{\rd\tau} \Big|_{\tau=\tau_*} t_{\alpha,\tau}\right) +1 \right) \right] e^{2 \pi i (t_{\alpha,\tau_*}+\tau_*)},
\end{equation}
 see equation \eqref{eq:formula-derivative-of-curve}. In particular, $\gamma_{\alpha,\rho}'(\tau_*)=0$ is equivalent to
\begin{equation}
\left[ \left( \frac{\rd}{\rd\tau} \Big|_{\tau=\tau_*} r_{\rho,\tau} \right)+ r_{\rho,\tau_*}  2\pi i\left(\left(\frac{\rd}{\rd\tau} \Big|_{\tau=\tau_*} t_{\alpha,\tau}\right) +1 \right) \right]=0.
\end{equation}
Therefore, from
\begin{equation}
\begin{aligned}
\gamma_{\alpha,\rho}''(\tau_*)&=\left[ \left( \frac{\rd^2}{\rd\tau^2} \Big|_{\tau=\tau_*} r_{\rho,\tau} \right)+ \left(\frac{\rd}{\rd\tau}\Big|_{\tau=\tau_*}r_{\rho,\tau}\right)  2\pi i\left(\left(\frac{\rd}{\rd\tau} \Big|_{\tau=\tau_*} t_{\alpha,\tau}\right) +1 \right)+\right.\\
& \qquad \left.+\,r_{\rho,\tau_*}  2\pi i\left(\frac{\rd^2}{\rd\tau^2} \Big|_{\tau=\tau_*} t_{\alpha,\tau}\right) \right] e^{2 \pi i (t_{\alpha,\tau_*}+\tau_*)}\\[1ex]
& \quad +\underbrace{ \left[ \left( \frac{\rd}{\rd\tau} \Big|_{\tau=\tau_*} r_{\rho,\tau} \right)+ r_{\rho,\tau_*}  2\pi i\left(\left(\frac{\rd}{\rd\tau} \Big|_{\tau=\tau_*} t_{\alpha,\tau}\right) +1 \right)\right]}_{=0}  \frac{\rd}{\rd\tau}\Big|_{\tau=\tau_*}e^{2 \pi i (t_{\alpha,\tau_*}+\tau_*)}\\[2ex]
&=\left[ \left( \frac{\rd^2}{\rd\tau^2} \Big|_{\tau=\tau_*} r_{\rho,\tau} \right)+ \left(\frac{\rd}{\rd\tau}\Big|_{\tau=\tau_*}r_{\rho,\tau}\right)  2\pi i\left(\left(\frac{\rd}{\rd\tau} \Big|_{\tau=\tau_*} t_{\alpha,\tau}\right) +1 \right)+\right.\\
& \qquad \left.+\,r_{\rho,\tau_*}  2\pi i\left(\frac{\rd^2}{\rd\tau^2} \Big|_{\tau=\tau_*} t_{\alpha,\tau}\right) \right] e^{2 \pi i (t_{\alpha,\tau_*}+\tau_*)}
\end{aligned}
\end{equation}
we conclude that $\gamma_{\alpha,\rho}''(\tau_*)\neq 0$ if and only if 
\begin{equation}
\frac{\rd^2}{\rd\tau^2} \Big|_{\tau=\tau_*} r_{\rho,\tau}\neq0\quad\text{or}\quad
\left(\frac{\rd}{\rd\tau}r_{\rho,\tau_*}\right) \left(\left(\frac{\rd}{\rd\tau} \Big|_{\tau=\tau_*} t_{\alpha,\tau}\right) +1 \right)+  r_{\rho,\tau_*} \left(\frac{\rd^2}{\rd\tau^2} \Big|_{\tau=\tau_*} t_{\alpha,\tau}\right)\neq0
\end{equation}
is non-zero. We claim that the first term never vanishes. Indeed, using \eqref{eq:def-of-r-rho-tau} we obtain
\begin{equation}\nonumber
\begin{aligned}
\frac{\rd^2}{\rd\tau^2} \Big|_{\tau=\tau_*} r_{\rho,\tau} &= \frac{\rd^2}{\rd\tau^2} \Big|_{\tau=\tau_*} \tilde{g}^{-1}_{\rho} \big(-2\pi \sin(2\pi\tau)\big) \\
&=\frac{\rd}{\rd\tau} \Big|_{\tau=\tau_*}\frac{-4 \pi^2 \cos(2\pi \tau_*)}{\tilde{g}'(r_{\rho,\tau_*})}\\
&=\frac{8\pi^3 \sin(2\pi \tau_*)}{\tilde{g}'(r_{\rho,\tau_*})} + \frac{4 \pi^2 \cos(2 \pi \tau_*)}{\tilde{g}'(r_{\rho,\tau_*})^2} \cdot \left(\tilde{g}''(r_{\rho,\tau_*} \right) \frac{\rd}{\rd \tau}\Big|_{\tau=\tau_*} r_{\rho,\tau}).
\end{aligned}
\end{equation}
Now, Lemma \ref{lem:app_cusp} asserts that $\gamma_{\alpha,\rho}'(\tau_*)=0$ implies $\tau_* \in \frac14+\Z$ or $\tau_*\in \frac34+\Z$. In either case the first summand in $\frac{\rd^2}{\rd\tau^2} \Big|_{\tau=\tau_*} r_{\rho,\tau}$ is non-zero and the second summand vanishes. This proves $\gamma_{\alpha,\rho}''(\tau_*)\neq 0$ and therefore the Proposition.
\end{proof}

With the analysis from the proof above at hand, we make some simple observations about the occurrence of cusps.

\begin{corollary}
    If the plot shows a cusp then it is always of 180°, meaning that the curve runs in and out of the cusp in opposite directions. 
\end{corollary}

\begin{corollary}
If $(z_{\alpha,\rho,\tau})_{\tau\in S^1}$ is a smooth $S^1$-family of delay orbits as in (2)(e) in Remark \ref{rem:collection_of_results_in_nondeg_case}   then there are at most two cusps in the curve $\gamma_{\alpha,\rho}$. Moreover, if there are two cusps then one of them satisfies condition (i) in Lemma~\ref{lem:app_cusp} but not (ii) and the other satisfies condition (ii) but not (i). 
\end{corollary}
 
\begin{proof}
We have seen in Lemma \ref{lem:app_cusp} that if cusps appear then $\tau_*\in\{\frac14,\frac34\} \subset S^1$ and condition (i) or (ii) in Lemma \ref{lem:app_cusp} need to hold. Recall from 2.(e) in Remark \ref{rem:collection_of_results_in_nondeg_case} that, by assumption, we have $f'(t_{\alpha,\tau})\neq0$ for all $\tau\in S^1$. If at both cusps condition (ii) in Lemma~\ref{lem:app_cusp} is satisfied then clearly $f'(t_{\alpha,\tau})=0$ for some $\tau\in S^1$, a contradiction. On the other hand, if we assume that condition (i) holds for both choices of $\tau_*$, that is, $r_{\rho,\frac{1}{4}}= r_{\rho,\frac{3}{4}} = 0$, then $\frac{\rd}{\rd\tau}r_{\rho,\tau}=\frac{-4 \pi^2 \cos(2\pi \tau)}{\tilde{g}'(r_{\rho,\tau})}=0$ for some $\tau$ in between $\frac14$ and $\frac34$ in $S^1$, another contradiction. See equation \eqref{eqn:formula_d_dtau_r_rho_tau} for the derivative of $r_{\rho,\tau}$. Therefore, for exactly one $\tau_*\in\{\frac14,\frac34\} \subset S^1$ satisfying condition (i) and the other satisfies condition (ii).
\end{proof}

\subsection{The degenerate case}
\label{subsec:deg_case}

Assume now that the 1-periodic orbit $z_{\alpha,\rho}$ (without delay) is degenerate, that is $f'(\alpha)=0$ for $\alpha\in Z=\{f=1\}$ or $g'(\rho)=0$ for $\rho\in R=\{g=0\}$. Note that the latter one implies $\tilde{g}'(\rho) = \frac{g'(\rho)\cdot\rho-g(\rho)}{\rho^2} = 0$. Similarly, $\tilde{g}$ has a maximum / minimum / saddle in $\rho>0$ if and only if $g$ does.

Next, following the construction of $t_{\alpha,\tau}\in S^1$ and $r_{\alpha,\rho}\in\R_{\geq 0}$ from Section~\ref{subsec:non-deg-case} and the ansatz $z_{\alpha,\rho,\tau}(t)=r_{\rho,\tau} e^{2 \pi i \cdot (t_{\alpha,\tau}+\tau+t)}$ we will point out similarities and difference in the degenerate cases.

\begin{itemize}
    \item If $\alpha$ or $\rho$ are saddle points of $f$ or $g$, then the construction of the orbits $z_{\alpha,\rho,\tau}$ works as in Lemma \ref{lem:delay-orbits}. The only difference from the non-degenerate case is that now the family $\{z_{\alpha,\rho,\tau}\}_{\tau\in[\tau_-,\tau_+]}$ 
    is not smooth for $\tau\in\Z$ if $f'(\alpha)=0$, see \eqref{eq:def-of-t-alpha-tau}, resp.~is not smooth in $\tau=\Z\cup(\frac12+\Z)$ if $\g'(\rho)=0$, see \eqref{eq:def-of-r-rho-tau}, since in these case the local inverses $f_\alpha^{-1}$ resp.~$\tilde{g}_\rho^{-1}$ are not smooth.
    \item If $f$ has a local minimum in $\alpha$ then there is no solution of \eqref{eq:equ-for-t-alpha-tau} for $\tau\not\in\Z$. Therefore we cannot find delay orbits using the ansatz $z_{\alpha,\rho,\tau}(t)=r_{\rho,\tau} e^{2 \pi i \cdot (t_{\alpha,\tau}+\tau+t)}$.
    \item If $f$ has a local maximum in $\alpha_{\text{max}}$, then there are two different choices for a local inverse of $f$ near $\alpha_{\text{max}}$. This yields one solution $t_{\alpha,\tau}$ for $\alpha$ smaller than (and close to) $\alpha_{\text{max}}$ and one solution for $\alpha$ larger than (and close to) $\alpha_{\text{max}}$. For both solutions, $\tau$ is close to $\Z$. In particular, we obtain two families of delay orbits meeting in $z_{\alpha,\rho}$ for $\tau\in\Z$. They are not smooth for $\tau\in\Z$.
    \item Similarly, if $f$ attains a local minimum in $\beta=t_{\alpha,\frac{1}{2}}\in S^1$ with $f(\beta)=-1$ then there are two different local inverses of $f$ near $-1$ with values around $\beta$. This implies that there are two families of delay orbits intersecting in the $\frac{1}{2}$-delay orbit $z_{\alpha,\rho,\frac{1}{2}}$. We are not aware of a suitable definition of non-degeneracy for delay orbits similar to Definition \ref{def:non_degenerate}, but this observation suggests that $z_{\alpha,\rho,\frac{1}{2}}$ should be called degenerate if $f'(t_{\alpha,\frac{1}{2}})=0$. Of course, here $\tfrac12$ may be replaced by any element in $\tfrac12+\Z$.
    \item Combining $\g(\rho)=0$ for  $\rho \in R$ with equation \eqref{eq:equ-for-r-rho-tau} we see that, if $\g$ attains a local minimum resp.~local maximum in $\rho \in R$, then our ansatz works only for $\tau<0$ resp. for $\tau>0$. 
    
    \item Similar to the case that $f$ attains a local minimum in $\beta=t_{\alpha,\frac{1}{2}}$ are the cases that $\tilde{g}$ attains a local minimum in $\lambda= r_{\rho,\frac{1}{4}}$ or a local maximum in $\kappa= r_{\rho,\frac{3}{4}}$. In these cases, there are two different choices for a local inverse of $\tilde{g}$ near $-2\pi$ with values around $\lambda$ or near $2\pi$ with values around $\kappa$ respectively. Hence there are two families of delay orbits meeting in the $\frac{1}{4}$-delay orbit $z_{\alpha,\rho,\frac{1}{4}}$ or in the $\frac{3}{4}$-delay orbit $z_{\alpha,\rho,\frac{3}{4}}$, respectively. Again, the orbits $z_{\alpha,\rho,\frac{1}{4}}$ and $z_{\alpha,\rho,\frac{3}{4}}$ should be called degenerate. Also, again we may add $\Z$ to $\tau=\frac14,\frac34$ in this discussion.
\end{itemize}

\section{Gluing families}

We have seen in Section~\ref{subsec:deg_case} that sometimes two families of delay orbits of the form $z_{\alpha,\rho,\tau}(t)=r_{\rho,\tau} e^{2 \pi i \cdot (t_{\alpha,\tau}+\tau+t)}$ meet in a degenerate orbit without delay, i.e. $\tau=0$ modulo $\Z$ or in a (degenerate) $\tau_*$-delay orbit with $\tau_* = \frac{1}{4}$, $ \frac{1}{2}$ or $\frac{3}{4}$ modulo $\Z$. More precisely:
\begin{itemize} \itemsep=1.5ex
\item $\tau_*=0\mod\Z$ if $f$ attains a local maximum in $\alpha=t_{\alpha,0}$ with $f(\alpha)=1$.
\item $\tau_*=\frac12\mod\Z$ if $f$ attains a local minimum in $\beta = t_{\alpha,\frac{1}{2}}$ with $f(\beta)=-1$.
\item $\tau_*=\frac14,\frac34\mod\Z$ if $\tilde{g}$ attains a local minimum in $\lambda=r_{\alpha,\frac{1}{4}}$ with $\tilde{g}(\lambda)=-2\pi$ or a maximum in $\kappa=r_{\alpha,\frac{3}{4}}$ with $\tilde{g}(\kappa)=2\pi$.
\end{itemize}
The two families then correspond to the two different choices of local inverse $f^{-1}_{\alpha,l}$ or $\tilde{g}^{-1}_{\kappa_j,l}$ ``to the left'' and $f^{-1}_{\alpha,r}$ or $\tilde{g}^{-1}_{\kappa_j,r}$ ``to the right'' of the respective extremum, which result in different solutions $t_{\alpha,\tau}$, $r_{\rho,\tau}$ to equations \eqref{eq:equ-for-t-alpha-tau} and  \eqref{eq:equ-for-r-rho-tau}.

In this section, we look at two specific classes of functions $f$ and $g$ for which explicit computations are possible and which lead to interesting examples. In Section \ref{subsec:trigonometric} we have a closer look at a particularly symmetric degenerate example coming from trigonometric functions where all extrema occur. We will see that it makes sense to glue the resulting families to each other in different ways and that we get very symmetric families of delay orbits from this glueing. After that, in Section \ref{subsec:polynomial}, we repeat the procedure for polynomial functions. The same glueing can of course be performed in examples with less symmetry.

\subsection[Trigonometric f and g]{Trigonometric $f$ and $g$}
\label{subsec:trigonometric}

We fix a positive integer $k \in \Z_{>0}$ and choose $f:S^1\to\R$ and $g:\R_{\geq0}\to\R$ as the following trigonometric functions.
\begin{align}
    f(\theta)&= \sin(2 \pi k \theta)   \\
    g(r) &= 2\pi r \cos( 2\pi r) \quad \text{i.e.} \quad \g(r)=2 \pi \cos(2 \pi r).
\end{align}
\begin{figure}[h]
    \centering
    \includegraphics[width=0.6\textwidth]{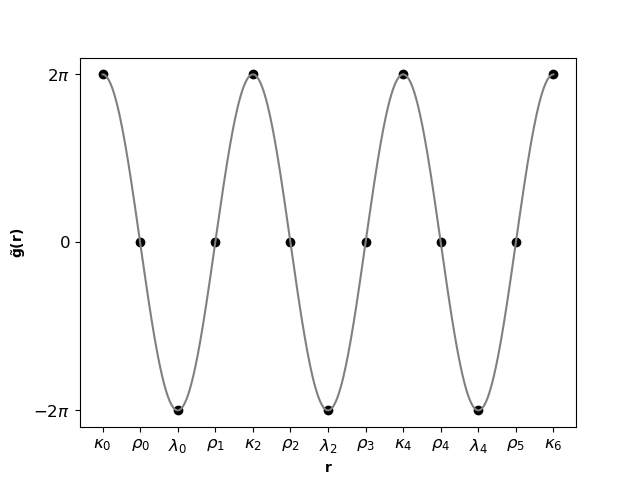}
    \caption{The function $\g$ with its zeros and extrema.}
    \label{fig:gluing_cos}
\end{figure}
In the notation from above we have $Z = f^{-1}(1) =\{\alpha_0,\dots,\alpha_{k-1}\}\subset S^1$ with $\alpha_j := \frac{4j+1}{4k}$,  $j=0,\ldots k-1$, and $f^{-1}(-1)=\{\beta_0,\dots,\beta_{k-1}\}$ with $\beta_j:=\frac{4j+3}{4k}$, $j=0,\ldots, k-1$. For $\g$ we have infinite sets of zeroes $R=\{0\} = \g^{-1}(0) = \{ \rho_j := \frac{2j+1}{4} \mid j\in \Z_{\geq 0}\}\subset\R$, maxima $\{\kappa_{\nu}:=\frac{\nu}{2} \mid \nu\in 2\Z_{\geq 0}\}$ and minima $\{\lambda_{\nu} := \frac{\nu+1}{2} \mid \nu\in 2\Z_{\geq 0}\}$, see Figure~\ref{fig:gluing_cos}. Note that we label maxima and minima by even numbers only. Then we have the following inequalities for $\nu\in 2\N_{\geq 0}$ between zeroes, maxima and minima of $\g$:
\begin{align*}
    \dots < \kappa_{\nu} < \rho_{\nu} < \lambda_{\nu} < \rho_{\nu+1} < \kappa_{\nu+2} < \rho_{\nu +2} < \lambda_{\nu+2} < \dots
\end{align*}
Near each maximum $\alpha_j\in Z$ of $f$, we have two choices for the local inverse of $f$. The inverse to the left of $\alpha_j=\frac{4j+1}{4k}$
\begin{align*}
        f_{\alpha_j,l}^{-1}:  [-1,1] &\longrightarrow [\beta_{j-1},\alpha_j]\subset S^1\\
        x & \longmapsto \frac{\arcsin(x)}{2 \pi k }+\frac{j}{k}
\end{align*}
and the inverse to the right of $\alpha_j$
\begin{align*}
        f_{\alpha_j,r}^{-1}: [-1,1] &\longrightarrow [\alpha_j,\beta_{j}]\subset S^1\\
        x & \longmapsto \frac{1}{2k}- \frac{\arcsin(x)}{2 \pi k }+\frac{j}{k}.
\end{align*}
For $\g$ we have near each zero $\rho_j= \frac{2j+1}{4}$ a local inverse which is, for even $j$, given by
\begin{align*}
    \g^{-1}_{\rho_j}: [-2\pi,2\pi] &\longrightarrow   [\kappa_j,\lambda_j] \\
    x &\longmapsto    \frac{\arccos(\frac{x}{2 \pi})}{2 \pi}+\frac{j}{2} 
\end{align*}
and for odd $j$ by
\begin{align*}
    \g^{-1}_{\rho_j}: [-2\pi,2\pi] &\longrightarrow  [\lambda_{j-1}, \kappa_{j+1}] \\
    x &\longmapsto      1-\frac{\arccos(\frac{x}{2 \pi})}{2 \pi}+\frac{j-1}{2}=-\frac{\arccos(\frac{x}{2 \pi})}{2 \pi}+\frac{j+1}{2}.
\end{align*}
Near the extrema of $\g$ we can choose between the local inverse $\g^{-1}_{\kappa_j,l}$ to the left and the local inverse $\g^{-1}_{\kappa_j,r}$ to the right. These are related to the other local inverse by $\g^{-1}_{\kappa_j,r} = \g^{-1}_{\rho_j} $ and $\g^{-1}_{\kappa_j,l}= \g^{-1}_{\rho_{j-1}}$ on the overlap of their respective domain of definition. To the left of a maximum $\kappa_{2j}$, $j\in \N_{\geq 0}$, we have the local inverse $\g^{-1}_{\rho_{2j-1}}$ (if $j>0$), to the right of $\kappa_{2j}$ we have $\g^{-1}_{\rho_{2j}}$. To the left of a minimum $\lambda_{2j}$, $j\in \N_{\geq 0}$, there is $\g^{-1}_{\rho_{2j}}$, to the right of $\lambda_{2j}$ there is $\g^{-1}_{\rho_{2j+1}}$.

All these local inverses are not smooth at the boundaries of the intervals of definition since their derivatives blow up. Therefore, if we use these inverses to define $S^1$-families of delay orbits via \eqref{eq:def-of-t-alpha-tau}, \eqref{eq:def-of-r-rho-tau} and \eqref{eq:def-of-delay-orbit} in Section \ref{subsec:non-deg-case}, we obtain $S^1$-families which are not smooth in $\tau=0,\frac{1}{2}$ due to $f_{\alpha_j,l}^{-1}$, and not smooth in $\tau=\frac{1}{4}, \frac{3}{4}$ due to $\g^{-1}_{\rho_j}$.

However, if we switch at each $\alpha_j$, $\beta_j$, $\lambda_j$ and $\kappa_j$ from the local inverse to the left to the local inverse to the right (or vice versa) as discussed at the beginning of this section we can define smooth families after all. Concretely, let us define solutions $t_{\alpha,\tau}^{lr}$ resp.~$t_{\alpha,\tau}^{rl}$ (${lr}$ and $rl$ indicates ``left to right'' and vice versa) to \eqref{eq:equ-for-t-alpha-tau} by
\begin{equation}\nonumber
\begin{aligned}
t^{lr}_{\alpha_j,\tau}&=
\begin{cases}
f_{\alpha_j,l}^{-1}(\cos(2\pi \tau)) & \tau \in[-\frac{1}{2},0] \\
f_{\alpha_j,r}^{-1}(\cos(2\pi \tau)) & \tau \in[0,\frac{1}{2}]
\end{cases} \\[1ex]
&=\frac{1}{k}(\tau +\frac{1}{4})+\frac{j}{k}
\end{aligned}
\end{equation}
respectively by
\begin{equation}\nonumber
\begin{aligned}
        t^{lr}_{\alpha_j,\tau}&=\begin{cases}
            f_{\alpha_j,r}^{-1}(\cos(2\pi \tau)) & \tau \in[-\frac{1}{2},0] \\
            f_{\alpha_j,l}^{-1}(\cos(2\pi \tau)) & \tau \in[0,\frac{1}{2}]
        \end{cases} \\[1ex]
        &=\frac{1}{2k}-\frac{1}{k}(\tau +\frac{1}{4})+\frac{j}{k}.
\end{aligned}
\end{equation}
Note that these families no longer form $S^1$-families but are parametrized over intervals. However, families for $j$ and $j+1$ agree in their ends, that is 
\begin{align*}
        t_{\alpha_j,\frac{1}{2}}^{lr} &=\beta_j = t_{\alpha_{j+1},-\frac{1}{2}}^{lr} \\
        t_{\alpha_j,-\frac{1}{2}}^{rl}    &=\beta_j = t_{\alpha_{j+1},\frac{1}{2}}^{rl}
\end{align*}
where $j=0,\ldots k-1$. Hence we obtain ``bigger'' families by glueing the ends of consecutive ``small'' families. In this way, we obtain two distinct families of solutions for equation \eqref{eq:equ-for-t-alpha-tau} given as
\begin{equation}\nonumber
   \begin{aligned}
        t^{lr}:\left[-\frac{1}{2},k-\frac{1}{2}\right] &\longrightarrow S^1 \\
     \tau & \longmapsto t_\tau^{lr} :=\frac{1}{k}\left(\tau + \frac{1}{4}\right) \\
    \end{aligned}\\
\end{equation}
and 
\begin{equation}\nonumber
    \begin{aligned}
        t^{rl}:\left[-\frac{1}{2},k-\frac{1}{2}\right] &\longrightarrow S^1 \\
     \tau & \longmapsto t_\tau^{rl}:= \frac{1}{k}\left(\frac{1}{4}-\tau \right). 
    \end{aligned}
\end{equation}
Both these solutions are $k$-periodic, so we can think of them as parameterized by $\tau \in \R/k\Z$. This reflects the fact that we started with a periodic function $f$.

Note that the corresponding two families of delay orbits do not depend on a particular choice of $\alpha\in Z\subset S^1$ anymore, instead, each family runs through every point in $Z$ consecutively. We arbitrarily chose to start at $\alpha_0=\frac{1}{4k}\in S^1 $ for $\tau=0$ for both solutions.

\begin{figure}[htb]
\captionsetup[subfigure]{justification=centering}
  \begin{subfigure}[t]{.45\textwidth}
    \centering
    \includegraphics[width=\linewidth]{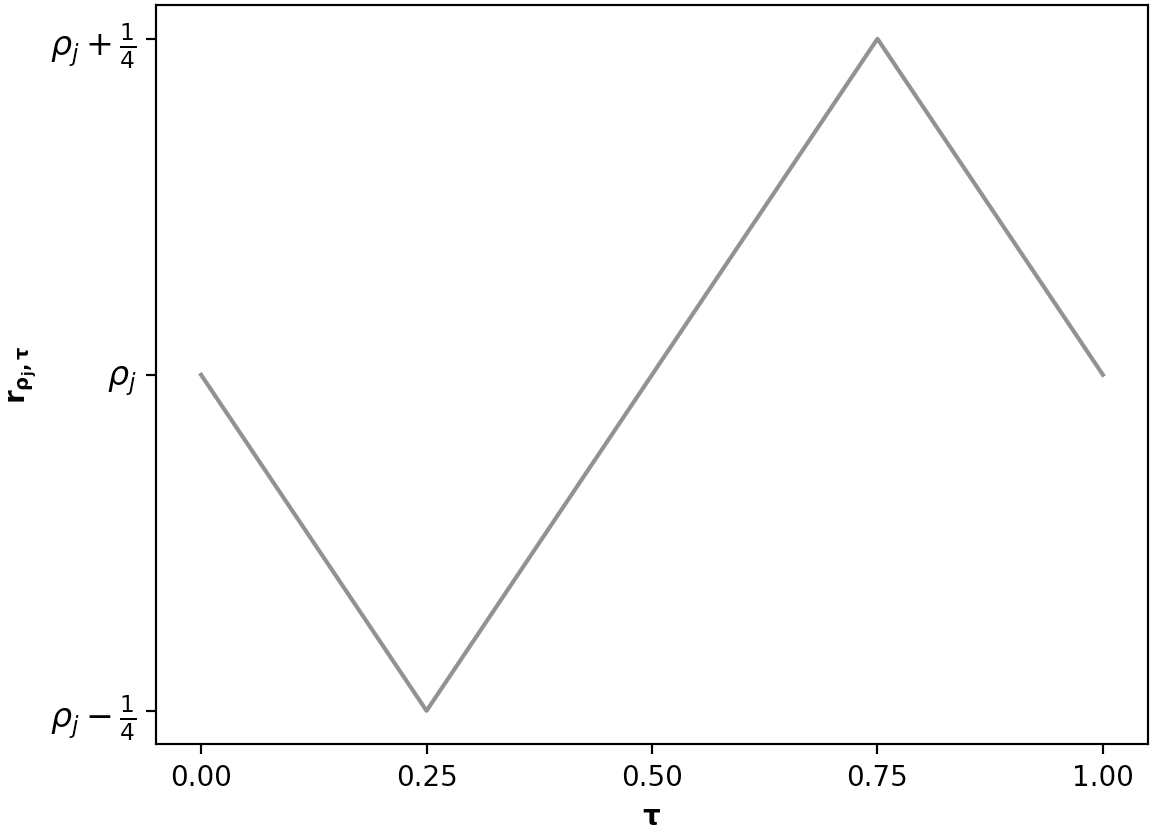}
    \caption{Mapping at $\rho_j$ for odd $j$}
    \label{fig:odd_inverse}
  \end{subfigure}
  \hfill
  \begin{subfigure}[t]{.45\textwidth}
    \centering
    \includegraphics[width=\linewidth]{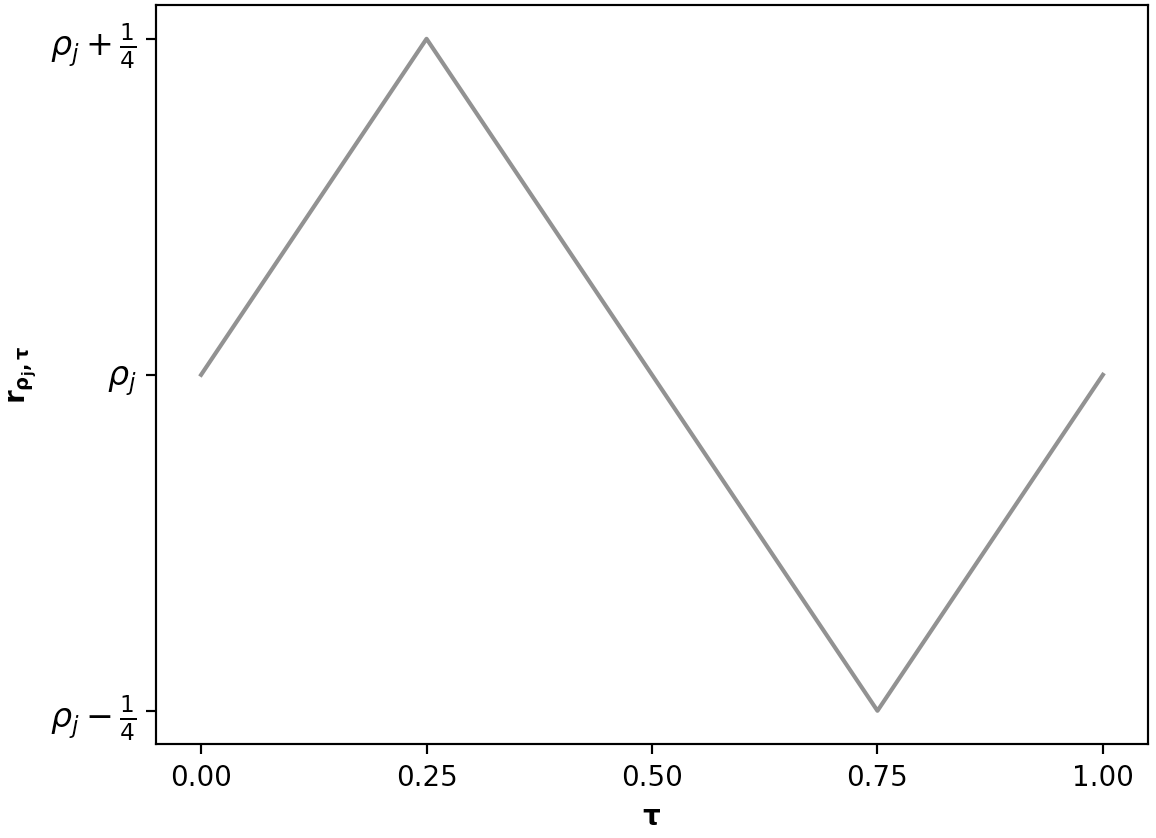}
    \caption{Mapping at $\rho_j$ for even $j$}
    \label{fig:even_inverse}
  \end{subfigure}
  \caption{Plot of $\tau \in [0,1] \rightarrow r_{\rho,\tau}$}
  \label{fig:gluing_radial_inverses}
\end{figure}

\begin{figure}[tb]
\captionsetup[subfigure]{justification=centering}
    \begin{subfigure}[t]{0.4\linewidth}
        \centering        \includegraphics[width=\linewidth,keepaspectratio]{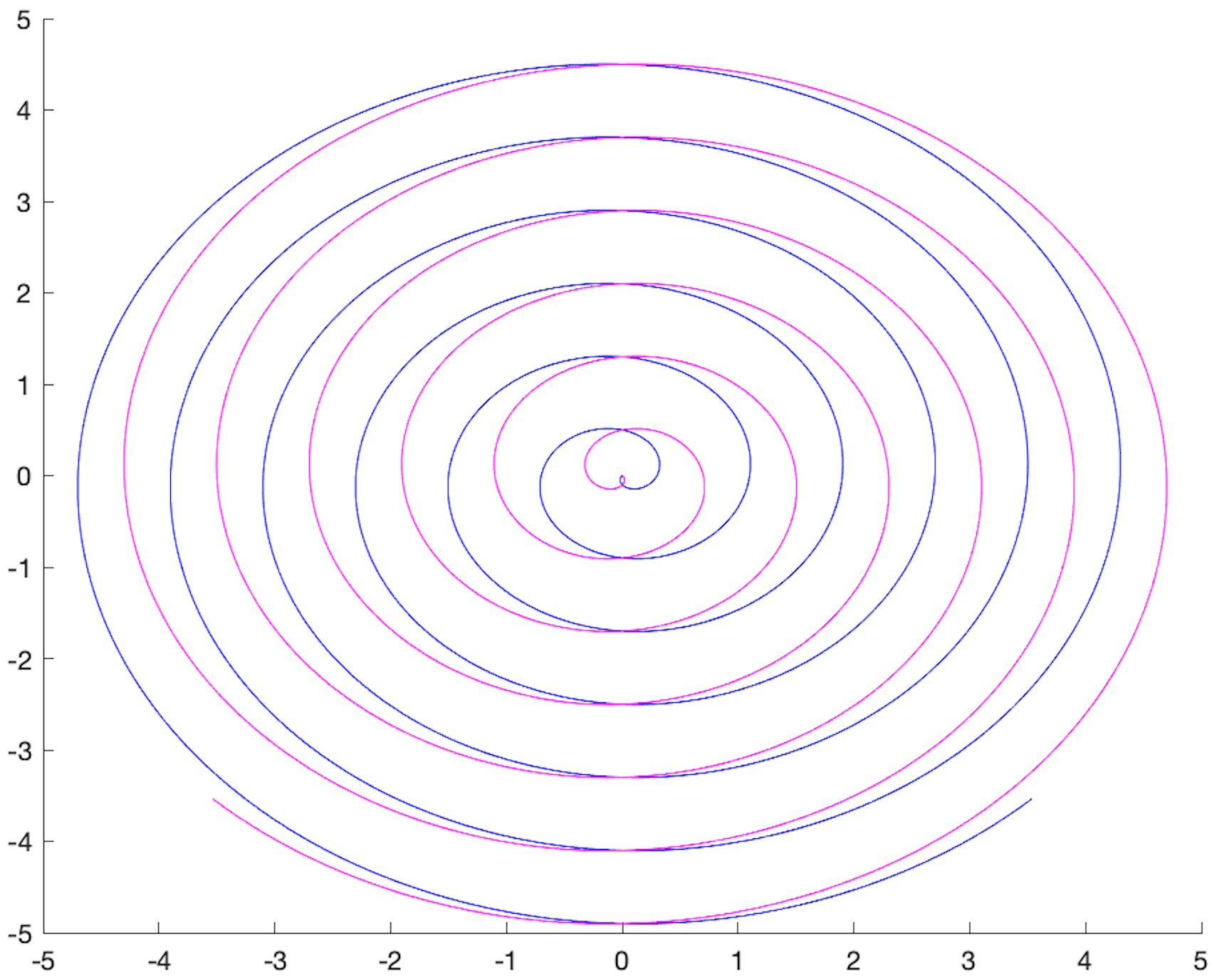}
         \caption{Combination of $r^{lr}$ (blue) \\ and $r^{rl}$ (magenta) with $t^{rl}$ }
         \label{fig:ex_tog_rl}
    \end{subfigure}
    \hfill
    \begin{subfigure}[t]{0.4\linewidth}
        \centering
         \includegraphics[width= \linewidth,keepaspectratio]{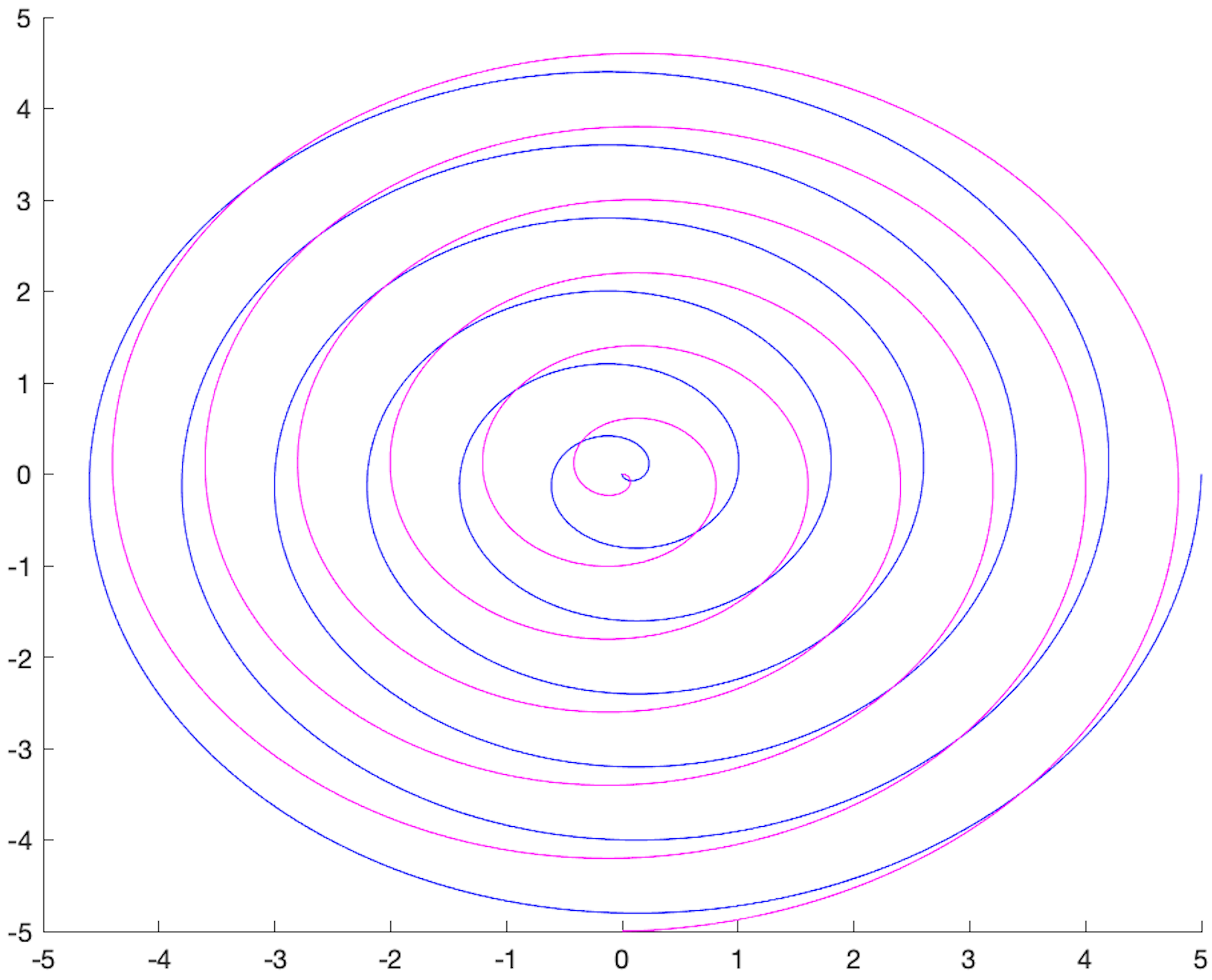}
         \caption{Combination of $r^{lr}$ (blue) \\ and $r^{rl}$ (magenta) with $t^{lr}$}
         \label{fig:ex_tog_lr}
    \end{subfigure}
    \caption{The four curves $\gamma: \tau \mapsto z_{\tau}^{lr/rl,\, lr/rl}(0) \in \R^2$ in the case of $k=5$, for $\tau \in[-\frac{1}{4},5+\frac{3}{4}]$ or $[-5-\frac{1}{4},\frac{3}{4}]$ respectively.}
\label{fig:gluing}
\end{figure}

Similarly, we proceed to define smooth families $r_{\tau}^{lr}$ and $r_{\tau}^{rl}$ of solutions to \eqref{eq:equ-for-r-rho-tau}. First, use the formula
\begin{align*}
    r_{\rho,\tau}:=\tilde{g}^{-1}_{\rho} \big(-2\pi \sin(2\pi\tau)\big)
\end{align*}
from \eqref{eq:def-of-r-rho-tau} for the different $\rho_j= \frac{2j+1}{4}$, $j\in \Z_{\geq 0}$. Since $\g(r)=2 \pi \cos(2 \pi r)$  this yields continuous, piecewise affine linear functions $\tau \mapsto r_{\rho_j,\tau}$ with $0\mapsto \rho_j$ and slope switching between $1$ and $-1$ at $\tau =\frac{1}{4},\frac{3}{4} \mod 1$, see Figure \ref{fig:gluing_radial_inverses}.
Then we can proceed and patch together the pieces for adjacent $j$. This would yield $\Z$-many affine linear (in particular, smooth) maps $\tau\mapsto r_{\rho_j,\tau}^{lr}$ (of slope $1$) resp.~$\tau \mapsto r_{\rho_j,\tau}^{rl}$ (of slope $-1$). However, the various $\tau\mapsto r_{\rho_j,\tau}^{lr}$ only differ by a shift from each other, stemming from the choice $r_{\rho_j,\tau=0}^{lr}=\rho_j$. Therefore, we prefer to define families $r_{\tau}^{lr}$ and $r_{\tau}^{rl}$ which do not depend on a particular point in $R$ but instead run through all points of $R$ consecutively again. We arbitrarily choose to let the ${lr}$-family pass through $\rho_0=\frac{1}{4}$ at $\tau=0$ and the $rl$-family through $\rho_1=\frac{3}{4}$ at $\tau=0$. This gives the following explicit expressions.
\begin{align*}
        r^{lr}:\R_{\geq \frac{-1}{4}} &\longrightarrow \R_{\geq 0} \\
        \tau &\longmapsto r_{\tau}^{lr} := \tau + \frac{1}{4} \\
        r^{rl}:\R_{\leq \frac{3}{4}} &\longrightarrow \R_{\geq 0} \\
        \tau &\longmapsto r_{\tau}^{lr} := \frac{3}{4} - \tau
\end{align*}
Finally, to define families of periodic delay orbits, we can now combine $t^{lr}$ and $t^{rl}$ with $r^{lr}$ and $r^{rl}$ resulting in four different families $z_{\tau}^{lr/rl,\, lr/rl}$, parametrized by $\tau\in \R_{\geq -\frac{1}{4}}$ or $\tau\in \R_{\leq \frac{3}{4}}$ respectively. Figure~\ref{fig:gluing} shows the four curves $\gamma: \tau \mapsto z_{\tau}^{lr/rl,\, lr/rl}(0) \in \R^2$ for $k=5$. For plotting reasons, we arbitrarily cut off the families of radii at $r_{\tau}^{lr,\,rl}=5$.

\subsection[Polynomial f and g]{Polynomial $f$ and $g$}
\label{subsec:polynomial}

\begin{figure}[h]
    \captionsetup[subfigure]{}%{justification=centering}
    \begin{subfigure}[t]{0.45\linewidth}
        \centering
         \includegraphics[width=\linewidth,keepaspectratio]{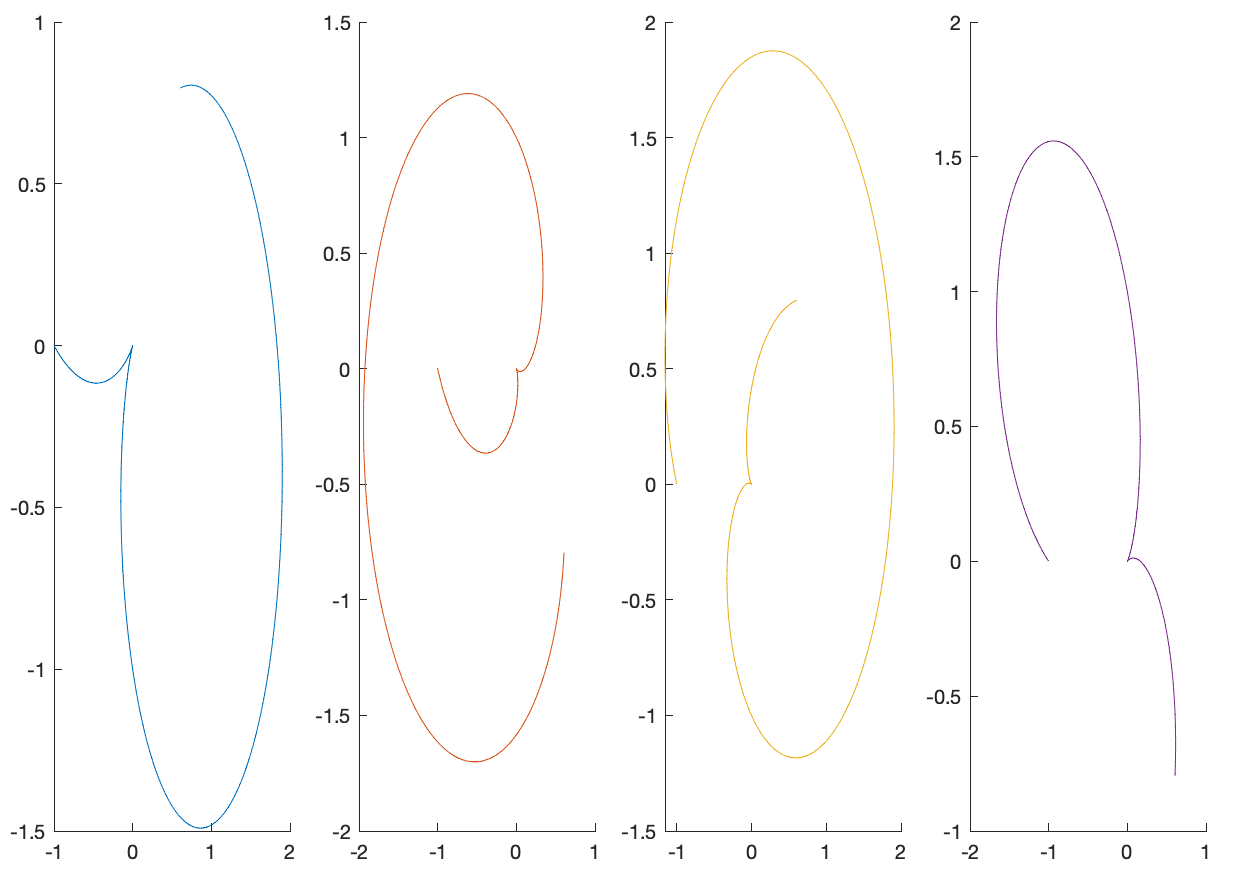}
         \caption{The four smooth non-closed curves $\gamma_{\alpha_i,\rho}$ having been glued only once.}
         \label{fig:glu_pol_df}
    \end{subfigure}
    \hfill
    \begin{subfigure}[t]{0.45\linewidth}
        \centering
        \includegraphics[width=\linewidth]{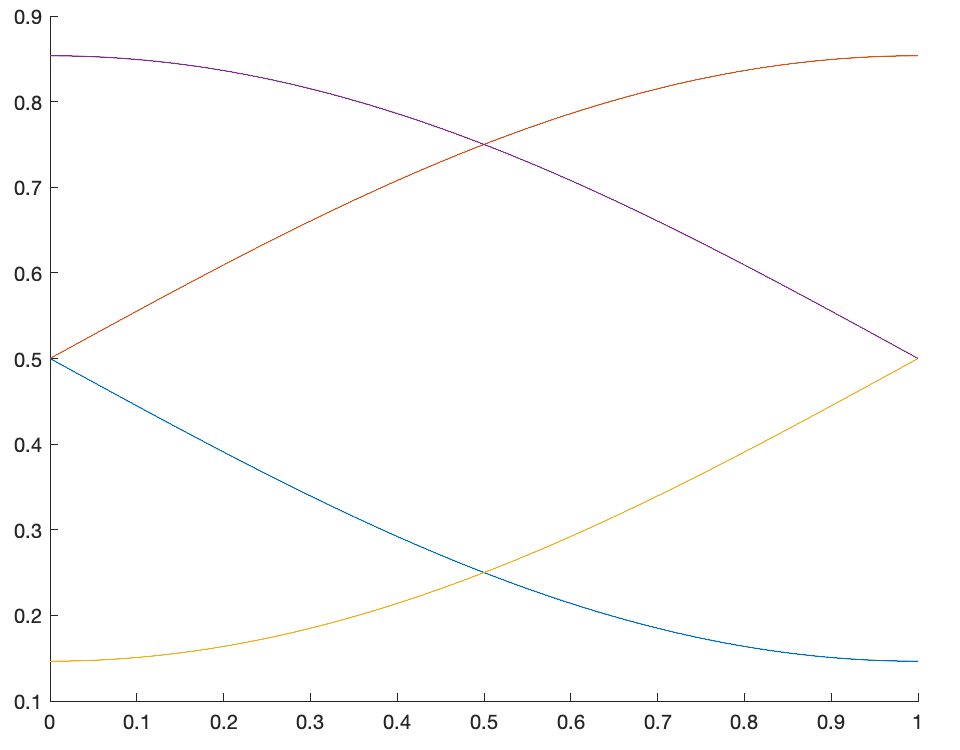}
        \caption{The glued $t_{\beta_j,\tau}$. The rising ones correspond to the solution $lr$ and the other way around.}
        \label{fig:glu_pol_t_al}
    \end{subfigure}
    \caption{Smooth delay families after first glueing.}
    \label{fig:glued_pol_1}
\end{figure}

\begin{figure}[h]
\captionsetup[subfigure]{justification=centering}
    \begin{subfigure}[t]{0.45\linewidth}
        \centering
        \includegraphics[width=\textwidth]{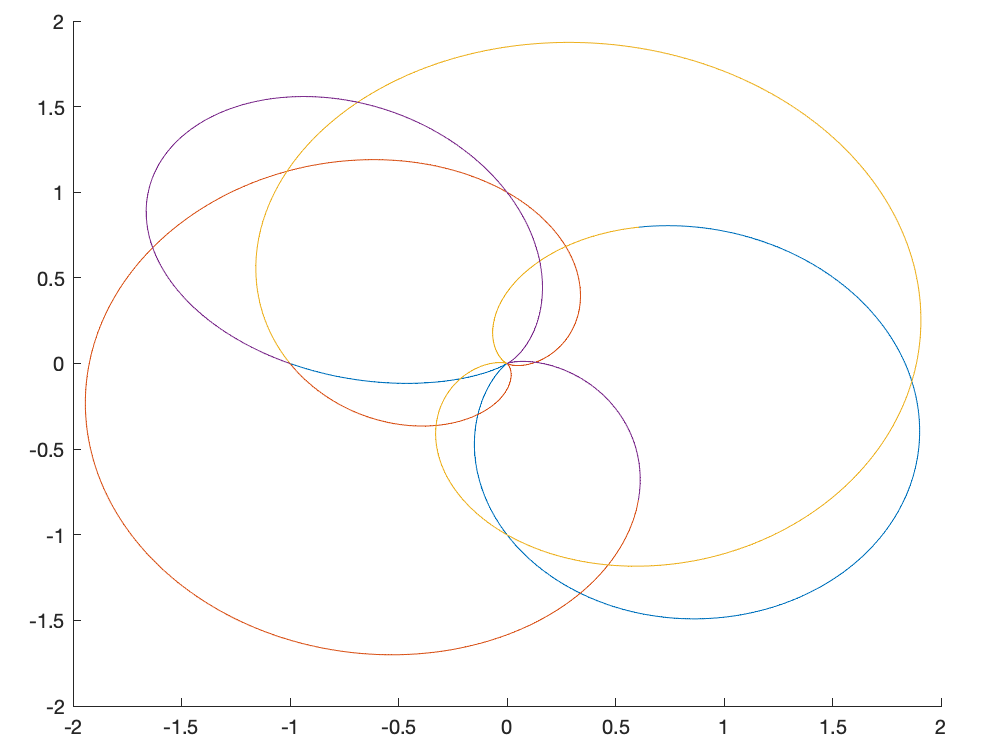}
        \caption{The glued delay family}
        \label{fig:glued_pol_2}
    \end{subfigure}
%    \hfill
    \begin{subfigure}[t]{0.45\linewidth}
        \centering
         \includegraphics[width= \linewidth,keepaspectratio]{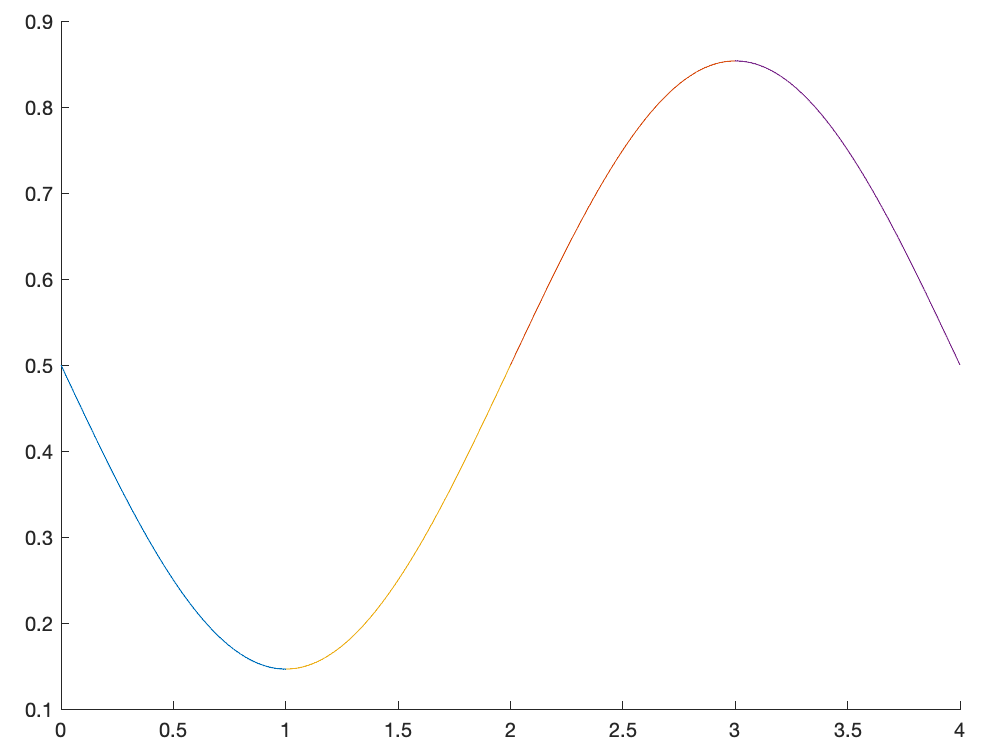}
         \caption{Plot of $[0,4]\to \R$, $\tau \mapsto t_{\text{glued}}(\tau) $}
         \label{fig:glued_ang_part}
    \end{subfigure}
    \caption{Same example as in Figure \ref{fig:glued_pol_1}, now all already glued families are glued together once more to form one smooth closed family.}
    \label{fig:tot_glued_angular_part}
\end{figure}

To showcase the above procedure in another example we consider
\begin{align*}
    f(\theta)&=512(\theta-\tfrac{1}{4})^2(\theta-\tfrac{3}{4})^2-1 \\
    g(r)&=2 \pi r(r-1).
\end{align*}
Clearly, $R=\{g=0\}=\{1\}$ and thus \eqref{eq:equ-for-r-rho-tau} leads to
\begin{align*}
    r_{\rho=1,\tau} = 1-\sin(2\pi \tau).
\end{align*}
From the expression for $f$, we see that we will have four distinct local inverses. These local inverses meet at the minima given by $\beta \in \{-\frac{1}{4},\frac{3}{4}\}=f^{-1}(-1)$ and the maximum at $\alpha \in \{\frac{1}{2}\} = f^{-1}(1)$. We denote the local inverses relative to their position to the two minima in $\beta_0=\frac{1}{4}$ resp.~$\beta_1 = \frac{3}{4}$ by
\begin{align*}
    f_{0,l}^{-1}&:[-1,1] \rightarrow [\alpha_0,\beta_0], \\
    f_{0,r}^{-1}&:[-1,1] \rightarrow [\beta_0, \alpha_1], \\
    f_{1,l}^{-1}&:[-1,1] \rightarrow [\alpha_1,\beta_1], \\
    f_{1,r}^{-1}&:[-1,1] \rightarrow [\beta_1, \alpha_2]. 
\end{align*}
Again, each of these local inverses on its own would not lead to a smoothly parameterized family of delay orbits. This can be seen from the kinks in Figure (\ref{fig:glu_pol_t_al}) which appear when switching between the red and violet arcs resp.~the yellow and blue arcs at $\tau=\frac12$.
However, this problem can be dealt with in the same manner as before, namely by combining the local inverses appropriately. We define the two left-to-right solutions by
\begin{align*}
        t^{lr}_{\beta_0,\tau}&=\begin{cases}
            f_{\beta_0,l}^{-1}(\cos(2\pi \tau)) & \tau \in[0,\frac{1}{2}] \\
            f_{\beta_0,r}^{-1}(\cos(2\pi \tau)) & \tau \in[\frac{1}{2},1]
        \end{cases} \\
        &=\frac{1}{2}- \frac{\sqrt{1+\cos(\pi \tau)}}{4}
\end{align*}
and
\begin{align*}
        t^{lr}_{\beta_1,\tau}&=\begin{cases}
            f_{\beta_1,l}^{-1}(\cos(2\pi \tau)) & \tau \in[0,\frac{1}{2}] \\
            f_{\beta_1,r}^{-1}(\cos(2\pi \tau)) & \tau \in[\frac{1}{2},1]
        \end{cases} \\
        &=\frac{1}{2}+ \frac{\sqrt{1-\cos(\pi \tau)}}{4}.
\end{align*}
The two right-to-left solutions are
\begin{align*}
        t^{rl}_{\beta_0,\tau}&=\begin{cases}
            f_{\beta_0,r}^{-1}(\cos(2\pi \tau)) & \tau \in[0,\frac{1}{2}] \\
            f_{\beta_0,l}^{-1}(\cos(2\pi \tau)) & \tau \in[\frac{1}{2},1]
        \end{cases} \\
        &=\frac{1}{2}-\frac{\sqrt{1-\cos(\pi \tau)}}{4}
\end{align*}
and
\begin{align*}
        t^{rl}_{\beta_1,\tau}&=\begin{cases}
            f_{\beta_1,r}^{-1}(\cos(2\pi \tau)) & \tau \in[0,\frac{1}{2}] \\
            f_{\beta_1,l}^{-1}(\cos(2\pi \tau)) & \tau \in[\frac{1}{2},1]
        \end{cases} \\
        &=\frac{1}{2}+ \frac{\sqrt{1+\cos(\pi \tau)}}{4}.
\end{align*}
Since $\cos(\pi) =-\cos(0)$, we can combine these into one smooth $S^1$-family
\begin{align*}
    t_{\text{glued}}(\tau)= \begin{cases}
        t_{\beta_0,\tau}^{rl} & \tau \in [0,1]\\
        t_{\beta_0,\tau}^{lr} & \tau \in [1,2]\\
        t_{\beta_1,\tau}^{lr} & \tau \in [2,3]\\
        t_{\beta_1,\tau}^{rl} & \tau \in [3,4].
    \end{cases}
\end{align*}
Any other cyclic permutation just gives another parameterization of the same family. The resulting $t_{\text{glued}}$ is shown in figure (\ref{fig:tot_glued_angular_part}).

\bibliographystyle{alpha}
\bibliography{bibliography.bib}

\end{document}